\documentclass[12pt]{article}

\usepackage{amsfonts,amssymb,amsthm,amsmath,latexsym}
\usepackage{subeqnarray,eqsection,indent,url,cite,bm}
\usepackage{tensor}  

\date{January 11, 2018 \\ revised July 3, 2018 \\[0.5mm]
      final version August 13, 2018 \\[1.5mm]
      To appear in {\em Expositiones Mathematicae}}   

\begin{document}

\title{The Euler and Springer numbers \\ as moment sequences}

\author{
     {\small Alan D.~Sokal}                  \\[2mm]
     {\small\it Department of Physics}       \\[-2mm]
     {\small\it New York University}         \\[-2mm]
     {\small\it 4 Washington Place}          \\[-2mm]
     {\small\it New York, NY 10003}          \\[-2mm]
     {\small\it USA}                         \\[-2mm]
     {\small\tt sokal@nyu.edu}               \\[-2mm]
     {\protect\makebox[5in]{\quad}}  
     \\[-2mm]
     {\small\it Department of Mathematics}   \\[-2mm]
     {\small\it University College London}   \\[-2mm]
     {\small\it Gower Street}                \\[-2mm]
     {\small\it London WC1E 6BT}             \\[-2mm]
     {\small\it UNITED KINGDOM}              \\[-2mm]
     {\small\tt sokal@math.ucl.ac.uk}        \\[3mm]
}

\maketitle
\thispagestyle{empty}   

\vspace*{-5mm}

\begin{abstract}
I study the sequences of Euler and Springer numbers
from the point of view of the classical moment problem.
\end{abstract}

\bigskip
\noindent
{\bf Key Words:}  Euler numbers, secant numbers, tangent numbers,
Springer numbers, alternating permutations, snakes of type $B_n$,
classical moment problem, Hamburger moment sequence, Stieltjes moment sequence,
Hankel matrix, Hankel determinant, continued fraction.

\bigskip
\bigskip
\noindent
{\bf Mathematics Subject Classification (MSC 2010) codes:}
05A15 (Primary);
05A05, 05A20, 06F25, 30B70, 30E05, 44A60, 60E99 (Secondary).

\clearpage

\newtheorem{theorem}{Theorem}[section]
\newtheorem{proposition}[theorem]{Proposition}
\newtheorem{lemma}[theorem]{Lemma}
\newtheorem{corollary}[theorem]{Corollary}
\newtheorem{definition}[theorem]{Definition}
\newtheorem{conjecture}[theorem]{Conjecture}
\newtheorem{question}[theorem]{Question}
\newtheorem{problem}[theorem]{Problem}
\newtheorem{example}[theorem]{Example}

\renewcommand{\theenumi}{\alph{enumi}}
\renewcommand{\labelenumi}{(\theenumi)}
\def\eop{\hbox{\kern1pt\vrule height6pt width4pt
depth1pt\kern1pt}\medskip}
\def\prf{\par\noindent{\bf Proof.\enspace}\rm}
\def\rmk{\par\medskip\noindent{\bf Remark\enspace}\rm}

\newcommand{\textbfit}[1]{\textbf{\textit{#1}}}

\newcommand{\bigdash}{%
\smallskip\begin{center} \rule{5cm}{0.1mm} \end{center}\smallskip}

\newcommand{\safepar}{ {\protect\hfill\protect\break\hspace*{5mm}} }

\newcommand{\be}{\begin{equation}}
\newcommand{\ee}{\end{equation}}
\newcommand{\<}{\langle}
\renewcommand{\>}{\rangle}
\newcommand{\widebar}{\overline}
\def\reff#1{(\protect\ref{#1})}
\def\spose#1{\hbox to 0pt{#1\hss}}
\def\ltapprox{\mathrel{\spose{\lower 3pt\hbox{$\mathchar"218$}}
    \raise 2.0pt\hbox{$\mathchar"13C$}}}
\def\gtapprox{\mathrel{\spose{\lower 3pt\hbox{$\mathchar"218$}}
    \raise 2.0pt\hbox{$\mathchar"13E$}}}
\def\textprime{${}^\prime$}
\def\proof{\par\medskip\noindent{\sc Proof.\ }}
\def\firstproof{\par\medskip\noindent{\sc First Proof.\ }}
\def\secondproof{\par\medskip\noindent{\sc Second Proof.\ }}
\def\alternateproof{\par\medskip\noindent{\sc Alternate Proof.\ }}
\def\algebraicproof{\par\medskip\noindent{\sc Algebraic Proof.\ }}
\def\combinatorialproof{\par\medskip\noindent{\sc Combinatorial Proof.\ }}
\def\proofof#1{\bigskip\noindent{\sc Proof of #1.\ }}
\def\firstproofof#1{\bigskip\noindent{\sc First Proof of #1.\ }}
\def\secondproofof#1{\bigskip\noindent{\sc Second Proof of #1.\ }}
\def\thirdproofof#1{\bigskip\noindent{\sc Third Proof of #1.\ }}
\def\algebraicproofof#1{\bigskip\noindent{\sc Algebraic Proof of #1.\ }}
\def\combinatorialproofof#1{\bigskip\noindent{\sc Combinatorial Proof of #1.\ }}
\def\sketchofproof{\par\medskip\noindent{\sc Sketch of proof.\ }}
\renewcommand{\qed}{ $\square$ \bigskip}
\newcommand{\myendremark}{ $\blacksquare$ \bigskip}
\def\half{ {1 \over 2} }
\def\third{ {1 \over 3} }
\def\twothird{ {2 \over 3} }
\def\smfrac#1#2{{\textstyle{#1\over #2}}}
\def\smhalf{ {\smfrac{1}{2}} }
\def\smquarter{ {\smfrac{1}{4}} }
\newcommand{\real}{\mathop{\rm Re}\nolimits}
\renewcommand{\Re}{\mathop{\rm Re}\nolimits}
\newcommand{\imag}{\mathop{\rm Im}\nolimits}
\renewcommand{\Im}{\mathop{\rm Im}\nolimits}
\newcommand{\sgn}{\mathop{\rm sgn}\nolimits}
\newcommand{\tr}{\mathop{\rm tr}\nolimits}
\newcommand{\tg}{\mathop{\rm tg}\nolimits}
\newcommand{\supp}{\mathop{\rm supp}\nolimits}
\newcommand{\disc}{\mathop{\rm disc}\nolimits}
\newcommand{\diag}{\mathop{\rm diag}\nolimits}
\newcommand{\csch}{\mathop{\rm csch}\nolimits}
\newcommand{\tridiag}{\mathop{\rm tridiag}\nolimits}
\newcommand{\AZ}{\mathop{\rm AZ}\nolimits}
\newcommand{\perm}{\mathop{\rm perm}\nolimits}
\def\hboxscript#1{ {\hbox{\scriptsize\em #1}} }
\renewcommand{\emptyset}{\varnothing}
\newcommand{\eqdef}{\stackrel{\rm def}{=}}

\newcommand{\restrict}{\upharpoonright}

\newcommand{\compinv}{{\langle -1 \rangle}}   

\newcommand{\scra}{{\mathcal{A}}}
\newcommand{\scrb}{{\mathcal{B}}}
\newcommand{\scrc}{{\mathcal{C}}}
\newcommand{\scrd}{{\mathcal{D}}}
\newcommand{\scre}{{\mathcal{E}}}
\newcommand{\scrf}{{\mathcal{F}}}
\newcommand{\scrg}{{\mathcal{G}}}
\newcommand{\scrh}{{\mathcal{H}}}
\newcommand{\scri}{{\mathcal{I}}}
\newcommand{\scrj}{{\mathcal{J}}}
\newcommand{\scrk}{{\mathcal{K}}}
\newcommand{\scrl}{{\mathcal{L}}}
\newcommand{\scrm}{{\mathcal{M}}}
\newcommand{\scrn}{{\mathcal{N}}}
\newcommand{\scro}{{\mathcal{O}}}
\newcommand{\scrp}{{\mathcal{P}}}
\newcommand{\scrq}{{\mathcal{Q}}}
\newcommand{\scrr}{{\mathcal{R}}}
\newcommand{\scrs}{{\mathcal{S}}}
\newcommand{\scrt}{{\mathcal{T}}}
\newcommand{\scrv}{{\mathcal{V}}}
\newcommand{\scrw}{{\mathcal{W}}}
\newcommand{\scrz}{{\mathcal{Z}}}

\newcommand{\bfa}{{\mathbf{a}}}
\newcommand{\bfb}{{\mathbf{b}}}
\newcommand{\bfc}{{\mathbf{c}}}
\newcommand{\bfd}{{\mathbf{d}}}
\newcommand{\bfe}{{\mathbf{e}}}
\newcommand{\bfj}{{\mathbf{j}}}
\newcommand{\bfi}{{\mathbf{i}}}
\newcommand{\bfk}{{\mathbf{k}}}
\newcommand{\bfl}{{\mathbf{l}}}
\newcommand{\bfm}{{\mathbf{m}}}
\newcommand{\bfx}{{\mathbf{x}}}
\renewcommand{\k}{{\mathbf{k}}}
\newcommand{\n}{{\mathbf{n}}}
\newcommand{\vv}{{\mathbf{v}}}
\newcommand{\bv}{{\mathbf{v}}}
\newcommand{\w}{{\mathbf{w}}}
\newcommand{\x}{{\mathbf{x}}}
\newcommand{\cc}{{\mathbf{c}}}
\newcommand{\zero}{{\mathbf{0}}}
\newcommand{\one}{{\mathbf{1}}}
\newcommand{\bmm}{{\mathbf{m}}}

\newcommand{\ahat}{{\widehat{a}}}
\newcommand{\Zhat}{{\widehat{Z}}}

\newcommand{\C}{{\mathbb C}}
\newcommand{\D}{{\mathbb D}}
\newcommand{\Z}{{\mathbb Z}}
\newcommand{\N}{{\mathbb N}}
\newcommand{\Q}{{\mathbb Q}}
\newcommand{\PP}{{\mathbb P}}
\newcommand{\R}{{\mathbb R}}
\newcommand{\RR}{{\mathbb R}}
\newcommand{\E}{{\mathbb E}}

\newcommand{\Sym}{{\mathfrak{S}}}
\newcommand{\SymB}{{\mathfrak{B}}}
\newcommand{\Alt}{{\mathrm{Alt}}}

\newcommand{\myle}{\preceq}
\newcommand{\myge}{\succeq}
\newcommand{\mygt}{\succ}

\newcommand{\B}{{\sf B}}
\newcommand{\OB}{{\sf OB}}
\newcommand{\OS}{{\sf OS}}
\newcommand{\OO}{{\sf O}}
\newcommand{\SP}{{\sf SP}}
\newcommand{\OSP}{{\sf OSP}}
\newcommand{\Eu}{{\sf Eu}}
\newcommand{\ERR}{{\sf ERR}}
\newcommand{\sfB}{{\sf B}}
\newcommand{\sfD}{{\sf D}}
\newcommand{\sfE}{{\sf E}}
\newcommand{\sfG}{{\sf G}}
\newcommand{\sfJ}{{\sf J}}
\newcommand{\sfP}{{\sf P}}
\newcommand{\sfQ}{{\sf Q}}
\newcommand{\sfS}{{\sf S}}
\newcommand{\sfT}{{\sf T}}
\newcommand{\sfW}{{\sf W}}
\newcommand{\sfMV}{{\sf MV}}
\newcommand{\AMV}{{\sf AMV}}
\newcommand{\BM}{{\sf BM}}
\newcommand{\NC}{{\sf NC}}

\newcommand{\emIB}{{\hbox{\em IB}}}
\newcommand{\emIP}{{\hbox{\em IP}}}
\newcommand{\emOB}{{\hbox{\em OB}}}
\newcommand{\emSC}{{\hbox{\em SC}}}

\newcommand{\stat}{{\rm stat}}
\newcommand{\cs}{{\rm cs}}
\newcommand{\cyc}{{\rm cyc}}
\newcommand{\Asc}{{\rm Asc}}
\newcommand{\asc}{{\rm asc}}
\newcommand{\Des}{{\rm Des}}
\newcommand{\des}{{\rm des}}
\newcommand{\Exc}{{\rm Exc}}
\newcommand{\exc}{{\rm exc}}
\newcommand{\Wex}{{\rm Wex}}
\newcommand{\wex}{{\rm wex}}
\newcommand{\Fix}{{\rm Fix}}
\newcommand{\fix}{{\rm fix}}
\newcommand{\lrmax}{{\rm lrmax}}
\newcommand{\rlmax}{{\rm rlmax}}
\newcommand{\Rec}{{\rm Rec}}
\newcommand{\rec}{{\rm rec}}
\newcommand{\Arec}{{\rm Arec}}
\newcommand{\arec}{{\rm arec}}
\newcommand{\ERec}{{\rm ERec}}
\newcommand{\erec}{{\rm erec}}
\newcommand{\EArec}{{\rm EArec}}
\newcommand{\earec}{{\rm earec}}
\newcommand{\recarec}{{\rm recarec}}
\newcommand{\nonrec}{{\rm nonrec}}
\newcommand{\Cpeak}{{\rm Cpeak}}
\newcommand{\cpeak}{{\rm cpeak}}
\newcommand{\Cval}{{\rm Cval}}
\newcommand{\cval}{{\rm cval}}
\newcommand{\Cdasc}{{\rm Cdasc}}
\newcommand{\cdasc}{{\rm cdasc}}
\newcommand{\Cddes}{{\rm Cddes}}
\newcommand{\cddes}{{\rm cddes}}
\newcommand{\cdrise}{{\rm cdrise}}
\newcommand{\cdfall}{{\rm cdfall}}
\newcommand{\Peak}{{\rm Peak}}
\newcommand{\peak}{{\rm peak}}
\newcommand{\Val}{{\rm Val}}
\newcommand{\val}{{\rm val}}
\newcommand{\Dasc}{{\rm Dasc}}
\newcommand{\dasc}{{\rm dasc}}
\newcommand{\Ddes}{{\rm Ddes}}
\newcommand{\ddes}{{\rm ddes}}
\newcommand{\inv}{{\rm inv}}
\newcommand{\maj}{{\rm maj}}
\newcommand{\rs}{{\rm rs}}
\newcommand{\cross}{{\rm cr}}
\newcommand{\crosshat}{{\widehat{\rm cr}}}
\newcommand{\nest}{{\rm ne}}
\newcommand{\rodd}{{\rm rodd}}
\newcommand{\reven}{{\rm reven}}
\newcommand{\lodd}{{\rm lodd}}
\newcommand{\leven}{{\rm leven}}
\newcommand{\sg}{{\rm sg}}
\newcommand{\bl}{{\rm bl}}
\newcommand{\tran}{{\rm tr}}
\newcommand{\area}{{\rm area}}
\newcommand{\ret}{{\rm ret}}
\newcommand{\peaks}{{\rm peaks}}
\newcommand{\hl}{{\rm hl}}
\newcommand{\sll}{{\rm sl}}
\newcommand{\negg}{{\rm neg}}
\newcommand{\imp}{{\rm imp}}
\newcommand{\osg}{{\rm osg}}
\newcommand{\ons}{{\rm ons}}
\newcommand{\isg}{{\rm isg}}
\newcommand{\ins}{{\rm ins}}
\newcommand{\LL}{{\rm LL}}
\newcommand{\height}{{\rm ht}}
\newcommand{\as}{{\rm as}}

\newcommand{\ba}{{\bm{a}}}
\newcommand{\bahat}{{\widehat{\bm{a}}}}
\newcommand{\sfa}{{{\sf a}}}
\newcommand{\bb}{{\bm{b}}}
\newcommand{\bc}{{\bm{c}}}
\newcommand{\bchat}{{\widehat{\bm{c}}}}
\newcommand{\bd}{{\bm{d}}}
\newcommand{\bee}{{\bm{e}}}
\newcommand{\bff}{{\bm{f}}}
\newcommand{\bg}{{\bm{g}}}
\newcommand{\bh}{{\bm{h}}}
\newcommand{\bll}{{\bm{\ell}}}
\newcommand{\bp}{{\bm{p}}}
\newcommand{\br}{{\bm{r}}}
\newcommand{\bs}{{\bm{s}}}
\newcommand{\bu}{{\bm{u}}}
\newcommand{\bw}{{\bm{w}}}
\newcommand{\bx}{{\bm{x}}}
\newcommand{\by}{{\bm{y}}}
\newcommand{\bz}{{\bm{z}}}
\newcommand{\bA}{{\bm{A}}}
\newcommand{\bB}{{\bm{B}}}
\newcommand{\bC}{{\bm{C}}}
\newcommand{\bE}{{\bm{E}}}
\newcommand{\bF}{{\bm{F}}}
\newcommand{\bG}{{\bm{G}}}
\newcommand{\bH}{{\bm{H}}}
\newcommand{\bI}{{\bm{I}}}
\newcommand{\bJ}{{\bm{J}}}
\newcommand{\bM}{{\bm{M}}}
\newcommand{\bN}{{\bm{N}}}
\newcommand{\bP}{{\bm{P}}}
\newcommand{\bQ}{{\bm{Q}}}
\newcommand{\bS}{{\bm{S}}}
\newcommand{\bT}{{\bm{T}}}
\newcommand{\bW}{{\bm{W}}}
\newcommand{\bX}{{\bm{X}}}
\newcommand{\bIB}{{\bm{IB}}}
\newcommand{\bOB}{{\bm{OB}}}
\newcommand{\bOS}{{\bm{OS}}}
\newcommand{\bERR}{{\bm{ERR}}}
\newcommand{\bSP}{{\bm{SP}}}
\newcommand{\bMV}{{\bm{MV}}}
\newcommand{\bBM}{{\bm{BM}}}
\newcommand{\balpha}{{\bm{\alpha}}}
\newcommand{\bbeta}{{\bm{\beta}}}
\newcommand{\bgamma}{{\bm{\gamma}}}
\newcommand{\bdelta}{{\bm{\delta}}}
\newcommand{\bkappa}{{\bm{\kappa}}}
\newcommand{\bomega}{{\bm{\omega}}}
\newcommand{\bsigma}{{\bm{\sigma}}}
\newcommand{\btau}{{\bm{\tau}}}
\newcommand{\bpsi}{{\bm{\psi}}}
\newcommand{\bzeta}{{\bm{\zeta}}}
\newcommand{\bone}{{\bm{1}}}
\newcommand{\bzero}{{\bm{0}}}

\newcommand{\Cbar}{{\overline{C}}}
\newcommand{\Dbar}{{\overline{D}}}
\newcommand{\dbar}{{\overline{d}}}
\def\Ctilde{{\widetilde{C}}}
\def\Etilde{{\widetilde{E}}}
\def\Ftilde{{\widetilde{F}}}
\def\Gtilde{{\widetilde{G}}}
\def\Htilde{{\widetilde{H}}}
\def\Ptilde{{\widetilde{P}}}
\def\Chat{{\widehat{C}}}
\def\ctilde{{\widetilde{c}}}
\def\zbar{{\overline{Z}}}
\def\pitilde{{\widetilde{\pi}}}

\newcommand{\sech}{{\rm sech}}

%
%
\newcommand{\sn}{{\rm sn}}
\newcommand{\cn}{{\rm cn}}
\newcommand{\dn}{{\rm dn}}
\newcommand{\sm}{{\rm sm}}
\newcommand{\cm}{{\rm cm}}

%
%
\newcommand{\zfz}{ {{}_0 \! F_0} }
\newcommand{\zfo}{ {{}_0  F_1} }
\newcommand{\ofz}{ {{}_1 \! F_0} }
\newcommand{\ofo}{ {{}_1 \! F_1} }
\newcommand{\oft}{ {{}_1 \! F_2} }

%
%
\newcommand{\FHyper}[2]{ {\tensor[_{#1 \!}]{F}{_{#2}}\!} }
\newcommand{\FHYPER}[5]{ {\FHyper{#1}{#2} \!\biggl(
   \!\!\begin{array}{c} #3 \\[1mm] #4 \end{array}\! \bigg|\, #5 \! \biggr)} }
\newcommand{\tfo}{ {\FHyper{2}{1}} }
\newcommand{\tfz}{ {\FHyper{2}{0}} }
\newcommand{\threefz}{ {\FHyper{3}{0}} }
\newcommand{\FHYPERbottomzero}[3]{ {\FHyper{#1}{0} \!\biggl(
   \!\!\begin{array}{c} #2 \\[1mm] \hbox{---} \end{array}\! \bigg|\, #3 \! \biggr)} }
\newcommand{\FHYPERtopzero}[3]{ {\FHyper{0}{#1} \!\biggl(
   \!\!\begin{array}{c} \hbox{---} \\[1mm] #2 \end{array}\! \bigg|\, #3 \! \biggr)} }

\newcommand{\phiHyper}[2]{ {\tensor[_{#1}]{\phi}{_{#2}}} }
\newcommand{\psiHyper}[2]{ {\tensor[_{#1}]{\psi}{_{#2}}} }
\newcommand{\PhiHyper}[2]{ {\tensor[_{#1}]{\Phi}{_{#2}}} }
\newcommand{\PsiHyper}[2]{ {\tensor[_{#1}]{\Psi}{_{#2}}} }
\newcommand{\phiHYPER}[6]{ {\phiHyper{#1}{#2} \!\left(
   \!\!\begin{array}{c} #3 \\ #4 \end{array}\! ;\, #5, \, #6 \! \right)\!} }
\newcommand{\psiHYPER}[6]{ {\psiHyper{#1}{#2} \!\left(
   \!\!\begin{array}{c} #3 \\ #4 \end{array}\! ;\, #5, \, #6 \! \right)} }
\newcommand{\PhiHYPER}[5]{ {\PhiHyper{#1}{#2} \!\left(
   \!\!\begin{array}{c} #3 \\ #4 \end{array}\! ;\, #5 \! \right)\!} }
\newcommand{\PsiHYPER}[5]{ {\PsiHyper{#1}{#2} \!\left(
   \!\!\begin{array}{c} #3 \\ #4 \end{array}\! ;\, #5 \! \right)\!} }
\newcommand{\zerophizero}{ {\phiHyper{0}{0}} }
\newcommand{\ophizero}{ {\phiHyper{1}{0}} }
\newcommand{\zphio}{ {\phiHyper{0}{1}} }
\newcommand{\ophio}{ {\phiHyper{1}{1}} }
\newcommand{\tphio}{ {\phiHyper{2}{1}} }
\newcommand{\tphiz}{ {\phiHyper{2}{0}} }
\newcommand{\tPhio}{ {\PhiHyper{2}{1}} }
\newcommand{\opsio}{ {\psiHyper{1}{1}} }

%
%
\newcommand{\stirlingsubset}[2]{\genfrac{\{}{\}}{0pt}{}{#1}{#2}}
\newcommand{\stirlingcycleold}[2]{\genfrac{[}{]}{0pt}{}{#1}{#2}}
\newcommand{\stirlingcycle}[2]{\left[\! \stirlingcycleold{#1}{#2} \!\right]}
\newcommand{\assocstirlingsubset}[3]{{\genfrac{\{}{\}}{0pt}{}{#1}{#2}}_{\! \ge #3}}
\newcommand{\genstirlingsubset}[4]{{\genfrac{\{}{\}}{0pt}{}{#1}{#2}}_{\! #3,#4}}
\newcommand{\euler}[2]{\genfrac{\langle}{\rangle}{0pt}{}{#1}{#2}}
\newcommand{\eulergen}[3]{{\genfrac{\langle}{\rangle}{0pt}{}{#1}{#2}}_{\! #3}}
\newcommand{\eulersecond}[2]{\left\langle\!\! \euler{#1}{#2} \!\!\right\rangle}
\newcommand{\eulersecondgen}[3]{{\left\langle\!\! \euler{#1}{#2} \!\!\right\rangle}_{\! #3}}
\newcommand{\binomvert}[2]{\genfrac{\vert}{\vert}{0pt}{}{#1}{#2}}
\newcommand{\binomsquare}[2]{\genfrac{[}{]}{0pt}{}{#1}{#2}}


\newenvironment{sarray}{
             \textfont0=\scriptfont0
             \scriptfont0=\scriptscriptfont0
             \textfont1=\scriptfont1
             \scriptfont1=\scriptscriptfont1
             \textfont2=\scriptfont2
             \scriptfont2=\scriptscriptfont2
             \textfont3=\scriptfont3
             \scriptfont3=\scriptscriptfont3
           \renewcommand{\arraystretch}{0.7}
           \begin{array}{l}}{\end{array}}

\newenvironment{scarray}{
             \textfont0=\scriptfont0
             \scriptfont0=\scriptscriptfont0
             \textfont1=\scriptfont1
             \scriptfont1=\scriptscriptfont1
             \textfont2=\scriptfont2
             \scriptfont2=\scriptscriptfont2
             \textfont3=\scriptfont3
             \scriptfont3=\scriptscriptfont3
           \renewcommand{\arraystretch}{0.7}
           \begin{array}{c}}{\end{array}}


\newcommand*\circled[1]{\tikz[baseline=(char.base)]{
  \node[shape=circle,draw,inner sep=1pt] (char) {#1};}}
\newcommand{\ostar}{{\circledast}}
\newcommand{\ostarN}{{\,\circledast_{\vphantom{\dot{N}}N}\,}}
\newcommand{\ostarPsi}{{\,\circledast_{\vphantom{\dot{\Psi}}\Psi}\,}}
\newcommand{\starN}{{\,\ast_{\vphantom{\dot{N}}N}\,}}
\newcommand{\starpsi}{{\,\ast_{\vphantom{\dot{\bpsi}}\!\bpsi}\,}}
\newcommand{\starone}{{\,\ast_{\vphantom{\dot{1}}1}\,}}
\newcommand{\startwo}{{\,\ast_{\vphantom{\dot{2}}2}\,}}
\newcommand{\starinfty}{{\,\ast_{\vphantom{\dot{\infty}}\infty}\,}}
\newcommand{\starT}{{\,\ast_{\vphantom{\dot{T}}T}\,}}

\clearpage

\section{Introduction}

The representation of combinatorial sequences as moment sequences
is a fascinating subject that lies at the interface between
combinatorics and analysis.  For instance, the Ap\'ery numbers
\be
   A_n
   \;=\;
   \sum_{k=0}^n \binom{n}{k}^{\! 2} \binom{n+k}{k}^{\! 2}
 \label{def.apery}
\ee
play a key role in Ap\'ery's celebrated proof \cite{Apery_79}
of the irrationality of $\zeta(3)$
\cite{vanderPoorten_79,Huylebrouck_01,Fischler_04,Zudilin_09};
they also arise in Ramanujan-like series for $1/\pi$
\cite{Chan_10}
and $1/\pi^2$ \cite{Zudilin_07}.
As such, they have elicited much interest,
both combinatorial \cite{Strehl_94,Schmidt_95,Chen_11b}
and number-theoretic \cite{Gessel_82,Beukers_85,Coster_90,Ahlgren_00,Osburn_16}.
A few years ago I conjectured \cite{Edgar_14},
based on extensive numerical computations,
that the Ap\'ery numbers are a Stieltjes moment sequence,
i.e.\ $A_n = \int \! x^n \, d\mu(x)$ for some positive measure $\mu$
on $[0,\infty)$.
Very recently this conjecture has been proven by Edgar \cite{Edgar_inprep},
in a {\em tour de force}\/ of special-functions work;
he gives an explicit formula, in terms of Heun functions,
for the (unique) representing measure~$\mu$.
The more general conjecture \cite{Sokal_totalpos} that the Ap\'ery polynomials
\be
   A_n(x)  \;=\;  \sum_{k=0}^n \binom{n+k}{k}^{\! 2} \binom{n}{k}^{\! 2} \, x^k
 \label{def.aperypoly}
\ee
are a Stieltjes moment sequence for all $x \ge 1$ remains open.

In this paper I propose to study the moment problem
for two less recondite combinatorial sequences:
the Euler numbers and the Springer numbers.
Many of the results given here are well known;
others are known but perhaps not as well known as they ought to be;
a few seem to be new.
This paper is intended as a leisurely survey
that presents the relevant results in a unified fashion
and employs methods that are as elementary as possible.

The \emph{Euler numbers} $E_n$
are defined by the exponential generating function
\be
   \sec t + \tan t  \;=\;  \sum_{n=0}^\infty E_n \, {t^n \over n!}
   \;.
 \label{def.euler}
\ee
The $E_{2n}$ are also called {\em secant numbers}\/,
and the $E_{2n+1}$ are called {\em tangent numbers}\/.
The Euler numbers are positive integers that satisfy the recurrence
\be
   E_{n+1}  \;=\;  {1 \over 2} \sum_{k=0}^n \binom{n}{k} E_{n-k} E_k
   \qquad\hbox{for $n \ge 1$}
 \label{eq.euler.recurrence}
\ee
with initial condition $E_0 = E_1 = 1$;
this recurrence follows easily from the differential equation
$\scre'(t) = \smhalf [1 + \scre(t)^2]$
for the generating function $\scre(t) = \sec t + \tan t$.
Andr\'e \cite{Andre_1879,Andre_1881} showed in 1879
that $E_n$ enumerates the alternating (down-up) permutations
of $[n] \eqdef \{1,\ldots,n\}$,
i.e.\ the permutations $\sigma \in \Sym_n$ that satisfy
$\sigma_1 > \sigma_2 < \sigma_3 > \sigma_4 < \ldots\;$.\footnote{
   As Josuat-Verg\`es {\em et al.}\/ point out
   \cite[p.~1613]{Josuat-Verges_12},
   Andr\'e's work ``is perhaps the first example of an inverse problem
   in the theory of generating functions:
   given a function whose Taylor series has nonnegative integer coefficients,
   find a family of combinatorial objects counted by those coefficients.''
}
More recently, other combinatorial objects have been found
to be enumerated by the Euler numbers:
complete increasing plane binary trees, increasing 0-1-2 trees,
Andr\'e permutations, simsun permutations, and many others;
see \cite{Foata_71,Viennot_81,Kuznetsov_94,Stanley_10} for surveys.
The sequence of Euler numbers starts as
\be
   (E_n)_{n \ge 0}
   \;=\;
   1, 1, 1, 2, 5, 16, 61, 272, 1385, 7936, 50521, \ldots
\ee
and can be found in the On-Line Encyclopedia of Integer Sequences \cite{OEIS}
as sequence~A000111.
By considering the singularities in the complex plane
of the generating function \reff{def.euler},
it is not difficult to show,
using a standard method \cite[pp.~258--259]{Flajolet_09},
that $E_n$ has the asymptotic behavior
\be
   E_n  \;=\;  {4 \over \pi} \, \Bigl( {2 \over \pi} \Bigr)^{\! n} \, n!
               \:+\: O\Bigl( \Bigl( {2 \over 3\pi} \Bigr)^{\! n} \, n! \Bigr)
 \label{eq.euler.asymptotics}
\ee
as $n \to\infty$.
See also \reff{eq.star3} below for a more precise
convergent expansion.\footnote{
   {\bf Warning:}
   The Euler numbers found in classical books of analysis
   are somewhat different from these:
   $E_{2n-1}^{\rm classical} = 0$
   and $E_{2n}^{\rm classical} = (-1)^n E_{2n}$.
   Moreover, our tangent numbers $E_{2n+1}$ are classically written
   as a complicated expression in terms of Bernoulli numbers:
$$
   E_{2n-1}  \;=\; {(-1)^{n-1} 2^{2n} (2^{2n} - 1) B_{2n}  \over 2n}
   \qquad\hbox{for $n \ge 1$}
   \;.
$$
   The definition given here is the one
   nowadays universally used by combinatorialists,
   since it makes $E_n$ a~positive integer that has a uniform combinatorial
   interpretation for $n$ even and $n$ odd.

   \quad
   Our Euler numbers can also be expressed in terms of the
   classical Euler polynomials $\scre_n(x)$
   defined by the generating function
$$
   \sum_{n=0}^\infty \scre_n(x) \, {t^n \over n!}
   \;=\;
   {2 e^{xt} \over e^t + 1}
   \;,
$$
   namely,
\begin{eqnarray*}
   E_{2n-1}  & = &  (-1)^n \, 2^{2n-1} \, \scre_{2n-1}(0)
             \;=\;  (-1)^{n-1} \, 2^{2n-1} \, \scre_{2n-1}(1)
             \qquad\hbox{for $n \ge 1$}  \\[2mm]
   E_{2n}    & = &  (-1)^n \, 2^{2n} \, \scre_{2n}(1/2)
\end{eqnarray*}
   I thank Christophe Vignat for this observation.

}

The {\em Springer numbers}\/ $S_n$
are defined by the exponential generating function
\cite{Glaisher_1898,Glaisher_14,Springer_71}
\be
   {1 \over \cos t - \sin t}
   \;=\;
   \sum_{n=0}^\infty S_n \, {t^n \over n!}
   \;.
 \label{def.springer}
\ee
Arnol'd \cite{Arnold_92} showed in 1992
that $S_n$ enumerates a signed-permutation analogue
of the alternating permutations.
More precisely, recall that a \emph{signed permutation} of $[n]$
is a sequence $\pi = (\pi_1,\ldots,\pi_n)$
of elements of $[\pm n] \eqdef \{-n,\ldots,-1\} \cup \{1,\ldots,n\}$
such that $|\pi| \eqdef (|\pi_1|,\ldots,|\pi_n|)$ is a permutation of $[n]$.
In~other words, a signed permutation $\pi$ is simply
a permutation $|\pi|$ together with a sign sequence $\sgn(\pi)$.
We write $\SymB_n$ for the set of signed permutations of $[n]$;
obviously $|\SymB_n| = 2^n n!$.
Then a {\em snake of type $B_n$}\/ is
a signed permutation $\pi \in \SymB_n$ that satisfies
$0 < \pi_1 > \pi_2 < \pi_3 > \pi_4 < \ldots\;$.
Arnol'd \cite{Arnold_92} showed 
that $S_n$ enumerates the snakes of type $B_n$.
Several other combinatorial objects are also
enumerated by the Springer numbers:
Weyl chambers in the principal Springer cone
of the Coxeter group $B_n$ \cite{Springer_71},
topological types of odd functions with $2n$ critical values \cite{Arnold_92},
and certain classes of complete binary trees and plane rooted forests
\cite{Josuat-Verges_14}.
The sequence of Springer numbers starts as
\be
   (S_n)_{n \ge 0}
   \;=\;
   1, 1, 3, 11, 57, 361, 2763, 24611, 250737, 2873041, 36581523, \ldots
\ee
and can be found in \cite{OEIS} as sequence~A001586.
It follows from \reff{def.springer} that $S_n$ has the asymptotic behavior
\be
   S_n  \;=\;  {2\sqrt{2} \over \pi} \, \Bigl( {4 \over \pi} \Bigr)^{\! n} \, n!
               \:+\: O\Bigl( \Bigl( {4 \over 3\pi} \Bigr)^{\! n} \, n! \Bigr)
 \label{eq.springer.asymptotics}
\ee
as $n \to\infty$.
See also \reff{eq.springer.star1} below for a more precise convergent expansion.

In this paper I propose to study the sequences of
Euler and Springer numbers from the point of view of
the classical moment problem
\cite{Stieltjes_1894,Shohat_43,Akhiezer_65,Simon_98,Schmudgen_17}.
Let us recall that
a sequence $\ba = (a_n)_{n \ge 0}$ of real numbers is called a
{\em Hamburger}\/ (resp.\ {\em Stieltjes}\/) {\em moment sequence}\/
if there exists a positive measure $\mu$
on $\R$ (resp.\ on $[0,\infty)$)
such that $a_n = \int \! x^n \, d\mu(x)$ for all $n \ge 0$.
A Hamburger (resp.\ Stieltjes) moment sequence is called
{\em H-determinate}\/ (resp.\ {\em S-determinate}\/)
if there is a unique such measure $\mu$;
otherwise it is called {\em H-indeterminate}\/ (resp.\ {\em S-indeterminate}\/).
Please note that a Stieltjes moment sequence can be
S-determinate but H-indeterminate
\cite[p.~240]{Akhiezer_65} \cite[p.~96]{Simon_98}.
The Hamburger and Stieltjes moment properties are also connected
with the representation of the ordinary generating function
$A(t) = \sum_{n=0}^\infty a_n t^n$
as a Jacobi-type or Stieltjes-type continued fraction;
this connection will be reviewed in Section~\ref{sec.prelim} below.

Many combinatorial sequences turn out to be Hamburger or Stieltjes
moment sequences,
and it is obviously of interest to find explicit expressions
for the representing measure(s) $\mu$
and/or the continued-fraction expansions of the ordinary generating function.
In this paper we will address both aspects for the Euler and Springer numbers
and some sequences related to them.

\section{Preliminaries on the moment problem}   \label{sec.prelim}

In this section we review some basic facts about the moment problem
\cite{Stieltjes_1894,Shohat_43,Wall_48,Akhiezer_65,Simon_98,Schmudgen_17}
that will be used repeatedly in the sequel.

In the Introduction we defined Hamburger and Stieltjes moment sequences.
We begin by noting some elementary consequences of these definitions:

1) If $\ba = (a_n)_{n \ge 0}$ is a Stieltjes moment sequence,
then every arithmetic-progression subsequence
$(a_{n_0 + jN})_{N \ge 0}$ with $n_0 \ge 0$ and $j \ge 1$
is again a Stieltjes moment sequence.

2) If $\ba = (a_n)_{n \ge 0}$ is a Hamburger moment sequence,
then every arithmetic-progression subsequence
$(a_{n_0 + jN})_{N \ge 0}$ with $n_0 \ge 0$ {\em even}\/ and $j \ge 1$
is again a Hamburger moment sequence;
and if also $j$ is even, then it is a Stieltjes moment sequence.

3) For a sequence $\ba = (a_n)_{n \ge 0}$, the following are equivalent:
\begin{quote}
\begin{itemize}
   \item[(a)]  $\ba$ is a Stieltjes moment sequence.
       \\[-5mm]
   \item[(b)]  The ``aerated'' sequence
       $\widehat{\ba} = (a_0,0,a_1,0,a_2,0,\ldots)$
       is a Hamburger moment sequence.
       \\[-5mm]
   \item[(c)]  There exist numbers $a'_0,a'_1,a'_2,\ldots$ such that
       the ``modified aerated'' sequence
       $\widehat{\ba}' = (a_0,a'_0,a_1,a'_1,a_2,a'_2,\ldots)$
       is a Hamburger moment sequence.
\end{itemize}
\end{quote}
Indeed, (b)$\implies$(c) is trivial,
and (c)$\implies$(a) follows from property \#2:
concretely, if $\widehat{\ba}'$ is represented by a measure $\widehat{\mu}'$
on $\R$,
then $\ba$ is represented by the measure $\mu$ on $[0,\infty)$
that is the image of $\widehat{\mu}'$ under the map $x \mapsto x^2$
[namely, $\mu(A) = \widehat{\mu}'(\{ x \colon\, x^2 \in A \})$].
And for (a)$\implies$(b),
if $\ba$ is represented by a measure $\mu$ supported on $[0,\infty)$,
then $\widehat{\ba}$ is represented by the even measure
$\widehat{\mu} = (\tau^+ + \tau^-)/2$ on $\R$,
where $\tau^\pm$ is the image of $\mu$
under the map $x \mapsto \pm \sqrt{x}$.

4) For a sequence $\ba = (a_n)_{n \ge 0}$, the following are equivalent:
\begin{quote}
\begin{itemize}
   \item[(a)]  $\ba$ is a Stieltjes moment sequence.
       \\[-5mm]
   \item[(b)]  Both $\ba$ and the once-shifted sequence
$\widetilde{\ba} = (a_{n+1})_{n \ge 0}$ are Stieltjes moment sequences.
       \\[-5mm]
   \item[(c)]  Both $\ba$ and $\widetilde{\ba}$ are Hamburger moment sequences.
\end{itemize}
\end{quote}
Here (a)$\iff$(b)$\implies$(c) is easy (using property \#1);
unfortunately I~do not know any completely elementary proof
of (c)$\implies$(a),
but it is anyway an immediate consequence of
Theorems~\ref{thm.hamburger} and \ref{thm.stieltjes} below
(see also \cite[p.~187]{Berg_84}).

5) If $\ba = (a_n)_{n \ge 0}$ and $\bb = (b_n)_{n \ge 0}$
are Hamburger (resp.\ Stieltjes) moment sequences,
then any linear combination $\alpha \ba + \beta \bb$
with $\alpha,\beta \ge 0$ is also a
Hamburger (resp.\ Stieltjes) moment sequence:
if $\ba$ (resp.\ $\bb$) has representing measure
$\mu$ (resp.\ $\nu$), then $\alpha \ba + \beta \bb$
has representing measure $\alpha\mu + \beta\nu$.

6) If $\ba = (a_n)_{n \ge 0}$ and $\bb = (b_n)_{n \ge 0}$
are Hamburger (resp.\ Stieltjes) moment sequences,
then their entrywise product $\ba\bb \eqdef (a_n b_n)_{n \ge 0}$
is also a Hamburger (resp.\ Stieltjes) moment sequence:
if $\ba$ (resp.\ $\bb$) has representing measure
$\mu$ (resp.\ $\nu$), then $\ba\bb$ has representing measure
given by the product convolution $\mu \diamond \nu$:
\be
   (\mu \diamond \nu)(A)
   \;=\;
   (\mu \times \nu)\bigl( \{(x,y) \in \R^2 \colon\: xy \in A \} \bigr)
   \qquad\hbox{for Borel $A \subseteq \R$}
\ee
[that is, $\mu \diamond \nu$ is the image of $\mu \times \nu$
 under the map $(x,y) \mapsto xy$].
We will often use this fact in the contrapositive:
if $\bb$ is a Hamburger (resp.\ Stieltjes) moment sequence
and $\ba\bb$ is {\em not}\/ a Hamburger (resp.\ Stieltjes) moment sequence,
then $\ba$ is {\em not}\/ a Hamburger (resp.\ Stieltjes) moment sequence.
Indeed, the non-Hamburger (resp.\ non-Stieltjes) property of $\ba\bb$
can be viewed as a strengthened form of the
non-Hamburger (resp.\ non-Stieltjes) property of $\ba$.

\bigskip

We now recall the well-known
\cite{Stieltjes_1894,Gantmakher_37,Shohat_43,Wall_48,Akhiezer_65,Simon_98,Schmudgen_17}
necessary and sufficient conditions for a sequence $\ba = (a_n)_{n \ge 0}$
to be a Hamburger or Stieltjes moment sequence.
To any infinite sequence $\ba = (a_n)_{n \ge 0}$ of real numbers,
we associate for each $m \ge 0$ the $m$-shifted infinite Hankel matrix
\be
   H_\infty^{(m)}(\ba)
   \;=\;
   (a_{i+j+m})_{i,j \ge 0}
   \;=\;
   \begin{bmatrix}
      a_m     & a_{m+1} & a_{m+2} & \cdots \\
      a_{m+1} & a_{m+2} & a_{m+3} & \cdots \\
      a_{m+2} & a_{m+3} & a_{m+4} & \cdots \\
      \vdots & \vdots  & \vdots & \ddots
   \end{bmatrix}
 \label{def.Hinftym.00}
\ee
and the $m$-shifted $n \times n$ Hankel matrix
\be
   H_n^{(m)}(\ba)
   \;=\;
   (a_{i+j+m})_{0 \le i,j \le n-1}
   \;=\;
   \begin{bmatrix}
      a_m & a_{m+1} & \cdots & a_{m+n-1}  \\
      a_{m+1} & a_{m+2} & \cdots & a_{m+n}  \\
      \vdots & \vdots  & \ddots & \vdots  \\
      a_{m+n-1} & a_{m+n} & \cdots & a_{m+2n-2}
   \end{bmatrix}
   \;.
 \label{def.Hnm.00}
\ee
We also define the Hankel determinants
\be
   \Delta_n^{(m)}(\ba) \;=\; \det H_n^{(m)}(\ba)
   \;.
\ee

\begin{theorem}[Necessary and sufficient conditions for Hamburger moment sequence]
   \label{thm.hamburger}
For a sequence $\ba = (a_n)_{n \ge 0}$ of real numbers,
the following are equivalent:
\begin{itemize}
   \item[(a)]  $\ba$ is a Hamburger moment sequence.
   \item[(b)]  $H_\infty^{(0)}(\ba)$ is positive-semidefinite.
      [That is, all the principal minors of $H_\infty^{(0)}(\ba)$
       are nonnegative.\footnote{
   Recall that a {\em minor}\/ of a matrix $A$
   is the determinant of a finite square submatrix $A_{IJ}$,
   where $I$ (resp.\ $J$) is a set of rows (resp.\ columns)
   and $|I| = |J| < \infty$.
   A {\em principal minor}\/ is the determinant of a finite submatrix $A_{II}$.
   See \cite[Theorem~7.2.5]{Horn_13} for the equivalence of
   positive-semidefiniteness with the nonnegativity of all principal minors.
}]
   \item[(c)]  There exist numbers
    $\alpha_0 \ge 0$, $\beta_1,\beta_2,\ldots \ge 0$ and
    $\gamma_0,\gamma_1,\ldots \in \R$ such that
\be
   \sum_{n=0}^{\infty} a_n t^n
   \;=\;
   \cfrac{\alpha_0}{1 - \gamma_0 t - \cfrac{\beta_1 t^2}{1 - \gamma_1 t - \cfrac{\beta_2 t^2}{1 - \cdots}}}
   \label{eq.thm.hankelreal.infty.Jtype.0}
\ee
in the sense of formal power series.
    [That is, the ordinary generating function
     $f(t) = \sum\limits_{n=0}^{\infty} a_n t^n$
    can be represented 
    as a Jacobi-type continued fraction with nonnegative coefficients $\bbeta$
    and $\alpha_0$.]
\end{itemize}
\end{theorem}

There is also a refinement that is often useful:
$\ba$ is a Hamburger moment sequence with a representing measure $\mu$
having infinite support
$\iff$
$H_\infty^{(0)}(\ba)$ is positive-definite
(i.e.\ all the principal minors are strictly positive)
$\iff$
all the leading principal minors $\Delta_n^{(0)}$ are strictly positive
$\iff$
all the $\beta_i$ are strictly positive.

\medskip

\begin{theorem}[Necessary and sufficient conditions for Stieltjes moment sequence]
   \label{thm.stieltjes}
For a sequence $\ba = (a_n)_{n \ge 0}$ of real numbers,
the following are equivalent:
\begin{itemize}
   \item[(a)]  $\ba$ is a Stieltjes moment sequence.
   \item[(b)]  Both $H^{(0)}_\infty(\ba)$ and $H^{(1)}_\infty(\ba)$
      are positive-semidefinite.
      [That is, all the principal minors of $H^{(0)}_\infty(\ba)$
       and $H^{(1)}_\infty(\ba)$ are nonnegative.]
   \item[(c)]  $H^{(0)}_\infty(\ba)$ is totally positive.
      [That is, {\em all} the minors of $H^{(0)}_\infty(\ba)$ are nonnegative.]
   \item[(d)]  There exist numbers $\alpha_0,\alpha_1,\ldots \ge 0$
      such that
\be
   \sum_{n=0}^{\infty} a_n t^n
   \;=\;
   \cfrac{\alpha_0}{1 - \cfrac{\alpha_1 t}{1 - \cfrac{\alpha_2 t}{1 - \cdots}}}
   \label{eq.thm.hankelreal.infty.Stype}
\ee
in the sense of formal power series.
    [That is, the ordinary generating function
     $f(t) = \sum\limits_{n=0}^{\infty} a_n t^n$
    can be represented 
    as a Stieltjes-type continued fraction with nonnegative coefficients.]
  \item[(e)]  There exist numbers
    $\alpha_0 \ge 0$, $\beta_1,\beta_2,\ldots \ge 0$ and
    $\gamma_0,\gamma_1,\ldots \ge 0$ such that
    the infinite tridiagonal matrix
\be
   A(\bbeta,\bgamma)
   \;=\;
   \begin{bmatrix}
      \gamma_0   & 1           &            &        &         \\
      \beta_1    & \gamma_1    & 1          &        &         \\
                 & \beta_2     & \gamma_2   & 1      &         \\
                 &             & \ddots     & \ddots & \ddots
   \end{bmatrix}
 \label{eq.tridiag}
\ee
is totally positive and
\be
   \sum_{n=0}^{\infty} a_n t^n
   \;=\;
   \cfrac{\alpha_0}{1 - \gamma_0 t - \cfrac{\beta_1 t^2}{1 - \gamma_1 t - \cfrac{\beta_2 t^2}{1 - \cdots}}}
   \label{eq.thm.hankelreal.infty.Jtype}
\ee
in the sense of formal power series.
    [That is, the ordinary generating function
     $f(t) = \sum\limits_{n=0}^{\infty} a_n t^n$
    can be represented 
    as a Jacobi-type continued fraction with a totally positive
    production matrix.]
\end{itemize}
\end{theorem}

Once again, there is a refinement:
$\ba$ is a Stieltjes moment sequence with a representing measure $\mu$
having infinite support
$\iff$
$H^{(0)}_\infty(\ba)$ and $H^{(1)}_\infty(\ba)$ are positive-definite
(i.e.\ all the principal minors are strictly positive)
$\iff$
all the leading principal minors $\Delta^{(0)}_n$ and $\Delta^{(1)}_n$
are strictly positive
$\iff$
$H^{(0)}_\infty(\ba)$ is strictly totally positive
(i.e.\ all the minors are strictly positive)
$\iff$
all the $\alpha_i$ are strictly positive
$\iff$
all the $\beta_i$ are strictly positive.

{}From the $2 \times 2$ minors of
$H^{(0)}_\infty(\ba)$ and $H^{(1)}_\infty(\ba)$,
we see that a Stieltjes moment sequence is log-convex:
$a_n a_{n+2} - a_{n+1}^2 \ge 0$.
(This is also easy to prove directly.)
But it goes without saying that
the Stieltjes moment property is {\em much}\/ stronger than log-convexity.

For future reference, let us also recall the formula
\cite[p.~21]{Wall_48} \cite[p.~V-31]{Viennot_83}
for the contraction of an S-fraction to a J-fraction:
\reff{eq.thm.hankelreal.infty.Stype} and \reff{eq.thm.hankelreal.infty.Jtype} 
are equal if
\begin{subeqnarray}
   \gamma_0  & = &  \alpha_1
       \slabel{eq.contraction_even.coeffs.a}   \\
   \gamma_n  & = &  \alpha_{2n} + \alpha_{2n+1}  \qquad\hbox{for $n \ge 1$}
       \slabel{eq.contraction_even.coeffs.b}   \\
   \beta_n  & = &  \alpha_{2n-1} \alpha_{2n}
       \slabel{eq.contraction_even.coeffs.c}
 \label{eq.contraction_even.coeffs}
\end{subeqnarray}

Concerning H-determinacy and S-determinacy,
we limit ourselves to quoting the following sufficient condition
\cite[Theorems~1.10 and 1.11]{Shohat_43}
due to Carleman in 1922:

\begin{theorem}[Sufficient condition for determinacy of moment problem]
   \label{thm.determinacy}
\quad\hfill
\vspace*{-2mm}
\begin{itemize}
   \item[(a)]  A Hamburger moment sequence $\ba = (a_n)_{n \ge 0}$
     satisfying $\sum\limits_{n=1}^\infty a_{2n}^{-1/2n} = \infty$
     is H-determinate.
   \item[(b)]  A Stieltjes moment sequence $\ba = (a_n)_{n \ge 0}$
     satisfying $\sum\limits_{n=1}^\infty a_{n}^{-1/2n} = \infty$
     is S-determinate.
\end{itemize}
\end{theorem}

\noindent
In Corollary~\ref{cor.determinacy} below,
we will prove, by elementary methods, a slightly weakened version
of Theorem~\ref{thm.determinacy}.
It should be stressed that the conditions of Theorem~\ref{thm.determinacy}
are sufficient for determinacy, but in no way necessary \cite{Heyde_63};
indeed, there are determinate Hamburger and Stieltjes moment sequences
with arbitrarily rapid growth \cite[pp.~89, 135]{Simon_98}.
In fact, 
given any H-indeterminate Hamburger (resp.\ Stieltjes) moment sequence
$\ba = (a_n)_{n \ge 0}$,
there exists an H-determinate Hamburger (resp.\ Stieltjes) moment sequence
$\ba' = (a'_n)_{n \ge 0}$
that differs from $\ba$ only in the zeroth entry:
$0 < a'_0 < a_0$ while $a'_n = a_n$ for all $n \ge 1$.\footnote{
   {\sc Proof} (for experts):
   If $\ba$ is an indeterminate Hamburger moment sequence,
   then the Nevanlinna-extremal measure corresponding to
   the parameter value $t=0$ (call it $\mu_0$)
   is a discrete measure concentrated on the zeros
   of the Nevanlinna $D$-function (which are all real and simple,
   and one of which is 0).
   If, in addition, $\ba$ is a Stieltjes moment sequence,
   then the orthonormal polynomials $P_n(x)$ have all their zeros
   in $(0,\infty)$, so $P_n(0) P_n(x) > 0$ for all $x \le 0$;
   it follows that $D(x) = x \sum_{n=0}^\infty P_n(0) P_n(x)$
   has all its zeros in $[0,\infty)$,
   so that $\mu_0$ is supported on $[0,\infty)$.
   Now consider the measure $\mu' = \mu_0 - \mu_0(\{0\}) \delta_0$:
   it is H-determinate \cite[p.~115]{Akhiezer_65} \cite[p.~111]{Berg_81}
   and its moment sequence $\ba'$ differs from $\ba$ only in the zeroth entry.
   I thank Christian Berg for drawing my attention to this result
   and its proof.
 \label{footnote_mu0}
}

We will need one other fact about determinacy \cite[p.~178]{Berg_94}:

\begin{proposition}[S-determinacy with H-indeterminacy]
   \label{prop.determinacy.berg}
Let $\ba$ be a Stieltjes moment sequence that is
S-determinate but H-indeterminate.
Then the unique measure on $[0,\infty)$ representing $\ba$
is the Nevanlinna-extremal measure corresponding to
the parameter value $t=0$,
hence is a discrete measure concentrated on the zeros
of the $D$-function from the Nevanlinna parametrization
(and in particular has an atom at 0).
\end{proposition}

\noindent
We refrain from explaining what is meant by ``Nevanlinna-extremal measure''
and ``Nevanlinna parametrization''
\cite{Akhiezer_65,Buchwalter_84b,Simon_98},
but simply stress that in this situation the representing measure
must be discrete.\footnote{
   {\sc Proof of Proposition~\ref{prop.determinacy.berg}} (for experts):
   Let $\ba$ be a Stieltjes moment sequence that is H-indeterminate.
   Then it was shown in footnote~\ref{footnote_mu0}
   that the N-extremal measure $\mu_0$
   is a discrete measure on $[0,\infty)$ representing $\ba$.
   If $\ba$ is also S-determinate, then $\mu_0$ is the
   unique measure on $[0,\infty)$ representing $\ba$.
   I again thank Christian Berg for drawing my attention to this result
   and its proof.
}

We will also make use of a generalization of the moment problem
from positive measures to signed measures.
So let $\mu$ be a finite signed measure on $\R$;
it has a unique Jordan decomposition $\mu = \mu_+ - \mu_-$
where $\mu_+, \mu_-$ are nonnegative and mutually singular \cite{Halmos_50}.
We write $|\mu| = \mu_+ + \mu_-$.
We will always assume that $|\mu|$ has finite moments of all orders,
i.e.\ $\int\limits_{-\infty}^\infty |x|^n \, d|\mu|(x) < \infty$
for all $n \ge 0$.
The moments $a_n = \int\limits_{-\infty}^\infty x^n \, d\mu(x)$
are then well-defined;
we say that $\mu$ represents $\ba = (a_n)_{n \ge 0}$.

In sharp contrast to Theorems~\ref{thm.hamburger} and \ref{thm.stieltjes},
the moment problem for signed measures has a trivial existence condition
and an extraordinary nonuniqueness:

\begin{theorem}[P\'olya \cite{Polya_37a,Polya_38}]
   \label{thm.polya}
Let $\ba = (a_n)_{n \ge 0}$ be \emph{any} sequence of real numbers,
and let $S$ be \emph{any} closed unbounded subset of $\R$.
Then there exists a signed measure $\mu$ with support in $S$
that represents $\ba$
[that is, $\int\limits_{-\infty}^\infty |x|^n \, d|\mu|(x) < \infty$
 for all $n \ge 0$
 and $a_n = \int\limits_{-\infty}^\infty x^n \, d\mu(x)$
 for all $n \ge 0$].
\end{theorem}

\noindent
So for any sequence $\ba$ (even the zero sequence!)\ there are
continuum many distinct signed measures $\mu$, with disjoint supports,
that represent $\ba$.
(For instance, we can take $S = \Z + \lambda$ for any $\lambda \in [0,1)$.)
See also Bloom \cite{Bloom_53} for a slight refinement;
and see Boas \cite{Boas_39} for a different proof of a weaker result.

The requirement here that $S$ be unbounded is essential;
among signed measures with bounded support, uniqueness holds.
More generally, uniqueness holds among signed measures that have
exponential decay.  To show this, we begin with some elementary lemmas:

\begin{lemma}[Bounded support]
   \label{lemma.bounded}
Let $\ba = (a_n)_{n \ge 0}$ be a sequence of real numbers,
let $\mu$ be a signed measure on $\R$ that represents $\ba$,
and let $R \in [0,\infty)$.
\begin{itemize}
   \item[(a)]  If $\mu$ is supported in $[-R,R]$,
      then $|a_n| \le \| \mu \| \, R^n$, where $\| \mu \| = |\mu|(\R)$.
   \item[(b)]  Conversely, if $\mu$ is a \emph{positive} measure
      and $|a_n| \le C R^n$ for some $C < \infty$,
      then $\mu$ is supported in $[-R,R]$.
\end{itemize}
\end{lemma}

\proof
(a) is trivial.

(b) Suppose that $\mu$ is a positive measure such that
$\mu \bigl( (-\infty,-R-\epsilon] \cup [R+\epsilon,\infty) \bigr) = K > 0$
for some $\epsilon > 0$.
Then $a_{2n} \ge K (R+\epsilon)^{2n}$ for all $n \ge 0$,
which contradicts the hypothesis $|a_n| \le C R^n$.
\qed

{\bf Remark.}
This proof shows that (b) holds under the weaker hypothesis
$\liminf\limits_{n \to\infty} |a_{2n}|^{1/2n} \le R$.
\myendremark

\begin{lemma}[Exponential decay]
   \label{lemma.exponential}
Let $\ba = (a_n)_{n \ge 0}$ be a sequence of real numbers,
let $\mu$ be a signed measure on $\R$ that represents $\ba$,
and let $\epsilon > 0$.
\begin{itemize}
   \item[(a)]  If
 $\displaystyle \int\limits_{-\infty}^\infty \! e^{\epsilon |x|} \, d|\mu|(x)
  = C < \infty$,
      then $|a_n| \le C \epsilon^{-n} n!$.
   \item[(b)]  Conversely, if $\mu$ is a \emph{positive} measure
      and $|a_n| \le C \epsilon^{-n} n!$ for some $C < \infty$,
      then
 ${\displaystyle \int\limits_{-\infty}^\infty \! e^{\delta |x|} \, d\mu(x)
  < \infty}$
      for all $\delta < \epsilon$.
\end{itemize}
\end{lemma}

\proof
(a) Since $|x^n| \le \epsilon^{-n} n! e^{\epsilon |x|}$,
it follows that $|a_n| \le C \epsilon^{-n} n!$.

(b) Applying the monotone convergence theorem
to $\cosh \delta x = \sum\limits_{n=0}^\infty (\delta x)^{2n}/(2n)!$,
we~conclude that
\be
   \int\limits_{-\infty}^\infty \! (\cosh \delta x) \, d\mu(x)
   \;=\;
   \sum\limits_{n=0}^\infty {\delta^{2n} \, a_{2n} \over (2n)!}
   \;\le\;
   {C \over 1 - \delta^2/\epsilon^2}
   \;<\;
   \infty
   \;.
\ee
\nopagebreak
\qed

\begin{proposition}[Uniqueness in the presence of exponential decay]
   \label{prop.uniqueness}
Let $\ba = (a_n)_{n \ge 0}$ be a sequence of real numbers,
and let $\mu$ and $\nu$ be signed measures on $\R$ that represent~$\ba$.
Suppose that $\mu$ has exponential decay in the sense that
$\displaystyle \int\limits_{-\infty}^\infty \! e^{\epsilon |x|} \, d|\mu|(x)
   < \infty$ for some $\epsilon > 0$;
and suppose that $\nu$ is either a positive measure
or else also has exponential decay.
Then $\mu = \nu$.
\end{proposition}

\proof
By Lemma~\ref{lemma.exponential}(a), we conclude that
$|a_n| \le C \epsilon^{-n} n!$ for some $C < \infty$.
Then Lemma~\ref{lemma.exponential}(b) implies that
if $\nu$ is a positive measure, it has exponential decay.
So we can assume that $\nu$ has exponential decay.
It follows that
$\displaystyle F(t) = \int\limits_{-\infty}^\infty \! e^{itx} \, d\mu(x)$
and
$\displaystyle G(t) = \int\limits_{-\infty}^\infty \! e^{itx} \, d\nu(x)$
define analytic functions in the strip
$|\!\imag t| < \epsilon$.
Moreover, by the dominated convergence theorem
they coincide in the disc $|t| < \epsilon$
with the absolutely convergent series
$\sum\limits_{n=0}^\infty a_n (it)^n/n!$.
It follows that $F=G$;
and by the uniqueness theorem for the Fourier transform
of tempered distributions \cite[Theorem~7.1.10]{Hormander_90}
(or by other arguments \cite[proof of Proposition~1.5]{Simon_98})
we conclude that $\mu = \nu$.
\qed

\begin{corollary}
   \label{cor.determinacy}
Let $\ba = (a_n)_{n \ge 0}$ be a sequence of real numbers
satisfying $|a_n| \le A B^n n!$ for some $A,B < \infty$.
Then there is at most one positive measure representing $\ba$.
\end{corollary}

\proof
Apply Lemma~\ref{lemma.exponential}(b) 
and then Proposition~\ref{prop.uniqueness}.
\qed

\begin{corollary}
   \label{cor.uniqueness}
Let $\ba = (a_n)_{n \ge 0}$ be a sequence of real numbers,
and let $\mu$ be a signed measure on $\R$
that is \emph{not} a positive measure,
that represents $\ba$, and that has exponential decay in the sense that
$\displaystyle \int\limits_{-\infty}^\infty \! e^{\epsilon |x|} \, d|\mu|(x)
   < \infty$ for some $\epsilon > 0$.
Then $\ba$ is \emph{not} a Hamburger moment sequence.
\end{corollary}

\section{Euler numbers, part 1}  \label{sec.euler1}

We begin by studying the sequence of Euler numbers divided by $n!$.
Our starting point is the partial-fraction expansions of secant and tangent
\cite[p.~11]{Andrews_99}:
\begin{eqnarray}
   \sec t
   & = &
  \lim\limits_{N \to \infty} \, \sum\limits_{k=-N}^N
        {(-1)^k \over (k+\smhalf)\pi - t}
     \label{eq.partfrac.sec}  \\[2mm]
   \tan t
   & = &
  \lim\limits_{N \to \infty} \, \sum\limits_{k=-N}^N
        {1 \over (k+\smhalf)\pi - t}
\end{eqnarray}
Inserting these formulae into the
exponential generating function \reff{def.euler} of the Euler numbers
and extracting coefficients of powers of $t$ on both sides, we obtain
\be
   {E_{2n} \over (2n)!}
   \;=\;
   \sum\limits_{k=-\infty}^\infty
   (-1)^k \, \bigl[ (k+\smhalf)\pi \bigr]^{-(2n+1)}
 \label{eq.star1}
\ee
(with the interpretation $\lim\limits_{N \to \infty} \, \sum\limits_{k=-N}^N$
 when $n=0$)
and
\be
   {E_{2n+1} \over (2n+1)!}
   \;=\;
   \sum\limits_{k=-\infty}^\infty
   \bigl[ (k+\smhalf)\pi \bigr]^{-(2n+2)}
   \;.
 \label{eq.star2}
\ee

We can rewrite \reff{eq.star2} as
\be
   {E_{2n+1} \over (2n+1)!}
   \;=\;
   \sum\limits_{k=0}^\infty
       {2 \over (k+\smhalf)^2 \pi^2} \:
       \biggl( {1 \over (k+\smhalf)^2 \pi^2} \biggr) ^{\! n}
   \;,
 \label{eq.star2a}
\ee
which represents $(E_{2n+1}/(2n+1)!)_{n \ge 0}$
as the moments of a positive measure supported on
a countably infinite subset of $[0, 4/\pi^2]$.
It follows that $(E_{2n+1}/(2n+1)!)_{n \ge 0}$ is a Stieltjes moment sequence,
which is both S-determinate and H-determinate.
Theorem~\ref{thm.stieltjes} then implies that
the ordinary generating function of $(E_{2n+1}/(2n+1)!)_{n \ge 0}$
can be written as a Stieltjes-type continued fraction
\reff{eq.thm.hankelreal.infty.Stype}
with nonnegative coefficients $\alpha_n$;
in fact we have the beautiful explicit formula \cite[p.~349]{Wall_48}
\be
   \sum_{n=0}^\infty {E_{2n+1} \over (2n+1)!} \, t^n
   \;=\;
   {\tan \sqrt{t} \over \sqrt{t}}
   \;=\;
   \cfrac{1}{1 - \cfrac{\frac{1}{3} t}{1 - \cfrac{\frac{1}{15} t}{1 -  \cfrac{\frac{1}{35} t}{1- \cdots}}}}
 \label{eq.lambert}
\ee
with coefficients $\alpha_n = 1/(4n^2 - 1) > 0$.
This continued-fraction expansion of the tangent function
was found by Lambert \cite{Lambert_1768} in 1761,
and used by him to prove the irrationality of $\pi$
\cite{Laczkovich_97,Wallisser_00}.
But in fact, as noted by Brezinski \cite[p.~110]{Brezinski_91},
a formula equivalent to \reff{eq.lambert} appears already
in Euler's first paper on continued fractions \cite{Euler_1744}:
see top p.~321 in the English translation.\footnote{
   The paper \cite{Euler_1744},
   which is E71 in Enestr\"om's \cite{Enestrom_13} catalogue,
   was presented to the St.~Petersburg Academy in 1737
   and published in 1744.
}
The expansion \reff{eq.lambert} is a $\zfo$ limiting case
of Gauss' continued fraction for the ratio of
two contiguous hypergeometric functions $\tfo$ \cite[Chapter~XVIII]{Wall_48}.

We will come back to \reff{eq.star1} in a moment.

Combining \reff{eq.star1} and \reff{eq.star2}
and taking advantage of the evenness/oddness of the summands, we get
\begin{subeqnarray}
   {E_n \over n!}
   & = &
   \sum\limits_{k=-\infty}^\infty
     \bigl[1 + (-1)^k \bigr] \, \bigl[(k+\smhalf)\pi \bigr]^{-(n+1)}
           \\[2mm]
   & = &
   2 \sum\limits_{k=-\infty}^\infty
       \biggl( {2 \over (4k+1)\pi} \biggr) ^{\! n+1}
 \slabel{eq.star3b}
 \label{eq.star3}
\end{subeqnarray}
(once again
 with the interpretation $\lim\limits_{N \to \infty} \, \sum\limits_{k=-N}^N$
 when $n=0$);
see \cite{Elkies_03} for further discussion of this sum.
For $n \ge 1$ this sum is absolutely convergent, so we can write
\be
   {E_{n+1} \over (n+1)!}
   \;=\;
   \sum\limits_{k=-\infty}^\infty
       {8 \over (4k+1)^2 \pi^2} \:
       \biggl( {2 \over (4k+1)\pi} \biggr) ^{\! n}
   \;,
 \label{eq.star3a}
\ee
which represents $(E_{n+1}/(n+1)!)_{n \ge 0}$
as the moments of a positive measure supported on
a countably infinite subset of $[-2/3\pi, 2/\pi]$.
It follows that $(E_{n+1}/(n+1)!)_{n \ge 0}$ is a Hamburger moment sequence,
which is H-determinate.
In fact, the ordinary generating function of $(E_{n+1}/(n+1)!)_{n \ge 0}$
can be written explicitly as a Jacobi-type continued fraction
\cite{Gladkovskii_12b}
\be
   \sum_{n=0}^\infty {E_{n+1} \over (n+1)!} \, t^n
   \;=\;
   \cfrac{1}{1 - \frac{1}{2}t - \cfrac{\frac{1}{12} t^2}{1 - \cfrac{\frac{1}{60} t^2}{1 -  \cfrac{\frac{1}{140} t^2}{1- \cdots}}}}
 \label{eq.Entilde1.Jfrac}
\ee
with coefficients $\gamma_0 = 1/2$, $\gamma_n = 0$ for $n \ge 1$,
and $\beta_n = 1/(16n^2 - 4)$.
This continued fraction can be obtained from
Lambert's continued fraction \reff{eq.lambert}
with $t$ replaced by $t^2/4$, by using the identity
\be
   \sum_{n=0}^\infty {E_{n+1} \over (n+1)!} \, t^n
   \;=\;
   {\sec t + \tan t - 1  \over  t}
   \;=\;
   {1 \over \displaystyle {t \over 2} \cot {t \over 2} \,-\, {t \over 2}}
   \;.
 \label{eq.Entilde1.identity}
\ee
By using the contraction formula \reff{eq.contraction_even.coeffs},
we can also rewrite \reff{eq.Entilde1.Jfrac}
as a Stieltjes-type continued fraction \cite{Gladkovskii_13b}
\be
   \sum_{n=0}^\infty {E_{n+1} \over (n+1)!} \, t^n
   \;=\;
   \cfrac{1}{1 - \cfrac{\frac{1}{2} t}{1 - \cfrac{\frac{1}{6} t}{1 +  \cfrac{\frac{1}{6} t}{1 + \cfrac{\frac{1}{10} t}{1- \cdots}}}}}
 \label{eq.Entilde1.Sfrac}
\ee
with coefficients $\alpha_{2k-1} = (-1)^{k-1}/(4k-2)$,
$\alpha_{2k} = (-1)^{k-1}/(4k+2)$.
Here the coefficients $\alpha_i$ are {\em not}\/ all nonnegative;
it follows by Theorem~\ref{thm.stieltjes}
and the uniqueness of Stieltjes-continued-fraction representations
that $(E_{n+1}/(n+1)!)_{n \ge 0}$
is {\em not}\/ a Stieltjes moment sequence
--- a fact that we already knew from \reff{eq.star3a}
and the H-determinacy.


Let us now consider the even subsequence $(E_{2n}/(2n)!)_{n \ge 0}$.
Is it a Hamburger moment sequence?
The answer is no, in a very strong sense:

\begin{proposition}
   \label{prop.E2nover2nfact}
Define $\Etilde_n = E_n/n!$.
Then $(\Etilde_{2n})_{n \ge 0}$ is not a Hamburger moment sequence.
In fact, no arithmetic-progression subsequence
$(\Etilde_{n_0 + jN})_{N \ge 0}$ with $n_0$ {\em even} and $j \ge 1$
is a Hamburger moment sequence.
\end{proposition}

\noindent
We give two proofs:

\firstproof
For any {\em even}\/ $n_0 \ge 2$ and any $j \ge 1$,
the equation~\reff{eq.star3a} represents
$(\widetilde{E}_{n_0 + jN})_{N \ge 0}$
as the moments of a signed measure supported on $[-2/3\pi, 2/\pi]$
that is {\em not}\/ a positive measure.
Corollary~\ref{cor.uniqueness} then implies that
$(\widetilde{E}_{n_0 + jN})_{N \ge 0}$
is not a Hamburger moment sequence.
The assertion for $n_0 = 0$ then follows from the assertion for $n_0 = 2j$.
\qed

The second proof is based on the following fact,
which is of some interest in its own right:

\begin{proposition}
   \label{prop.euler.polya}
Define $\Etilde_n = E_n/n!$.
Then $(\Etilde_{2n})_{n \ge 0}$ is a P\'olya frequency sequence,
i.e.\ every minor of the infinite Toeplitz matrix
\be
   (\Etilde_{2j-2i})_{i,j \ge 0}
   \;=\;
   \begin{bmatrix}
      \Etilde_0 & \Etilde_2 & \Etilde_4 & \Etilde_6 & \cdots \\
      0         & \Etilde_0 & \Etilde_2 & \Etilde_4 & \cdots \\
      0         & 0         & \Etilde_0 & \Etilde_2 & \cdots \\
      0         & 0         & 0         & \Etilde_0 & \cdots \\
      \vdots    & \vdots    & \vdots    & \vdots    & \ddots
   \end{bmatrix}
\ee
is nonnegative.
Moreover, a minor using rows $i_1 < i_2 < \ldots < i_r$
and columns $j_1 < j_2 < \ldots < j_r$
is strictly positive if $i_k \le j_k$ for $1 \le k \le r$.
In particular, the sequence $(\Etilde_{2n})_{n \ge 0}$
is strictly log-concave.
\end{proposition}

\proof
It follows from the well-known infinite product representation
for $\cos t$ that
\be
   \sum_{n=0}^\infty \Etilde_{2n} \, t^n
   \;=\;
   \sec \sqrt{t}
   \;=\;
   \prod_{k=0}^\infty
       \biggl( 1 \,-\, {t \over (k+\smhalf)^2 \pi^2} \biggr) ^{\! -1}
   \;.
\ee
This implies \cite[p.~395]{Karlin_68}
that $(\Etilde_{2n})_{n \ge 0}$ is a P\'olya frequency sequence;
and it also implies \cite[p.~427--430]{Karlin_68}
the statement about strictly positive minors.
The strict log-concavity is simply the strict positivity
of the $2 \times 2$ minors above the diagonal.
\qed

\vspace*{-4mm}

\secondproofof{Proposition~\ref{prop.E2nover2nfact}}
No arithmetic-progression subsequence
$(\widetilde{E}_{n_0 + jN})_{N \ge 0}$ with $n_0$ even and $j \ge 1$
can be a Hamburger moment sequence,
since its even subsequence $(\Etilde_{n_0 +2jN})_{N \ge 0}$
is strictly log-concave and hence cannot be log-convex.
\qed

The even and odd subsequences thus have radically different behavior:
the even subsequence $(\Etilde_{2n})_{n \ge 0}$ is strictly log-concave
(since it is a P\'olya frequency sequence with an infinite representing
 product),
while the odd subsequence $(\Etilde_{2n+1})_{n \ge 0}$ is strictly log-convex
(since it is a Stieltjes moment sequence with a representing measure
of infinite support).
These two facts are special cases of the following more general inequality
that appears to be true:

\begin{conjecture}
   \label{conj.euler.logconv}
Define $\Etilde_n = E_n/n!$.
Then for all $n \ge 0$ and $j,k \ge 1$, we have
\be
   (-1)^{n-1} [\Etilde_n \Etilde_{n+j+k} - \Etilde_{n+j} \Etilde_{n+k}]
   \;>\;  0 
   \;.
 \label{eq.conj.euler.logconv}
\ee
\end{conjecture}

\noindent
I do not know how to prove \reff{eq.conj.euler.logconv},
but I have verified it for $n,j,k \le 900$.



%

Though $(E_{2n}/(2n)!)_{n \ge 0}$ is not a Hamburger moment sequence,
one could try multiplying it by a Hamburger (or Stieltjes) moment sequence
$(b_n)_{n \ge 0}$;
the result $(b_n E_{2n}/(2n)!)_{n \ge 0}$ {\em might}\/
be a Hamburger (or even a Stieltjes) moment sequence.
For instance, the central binomial coefficients
$\binom{2n}{n} = (2n)!/(n!)^2$ are a Stieltjes moment sequence,
with representation
\be
   \binom{2n}{n}
   \;=\;
   {1 \over \pi} \int\limits_0^4 x^n \, x^{-1/2} \, (4-x)^{-1/2} \, dx
 \label{eq.centralbinomial.stieltjes}
\ee
as a special case of the beta integral.
Might $(E_{2n}/(n!)^2)_{n \ge 0}$ be a Hamburger moment sequence?
The answer is no, because the $7 \times 7$ Hankel determinant
$\det (a_{i+j})_{0 \le i,j \le 6}$ for $a_n = E_{2n}/(n!)^2$ is negative.
Unfortunately I do not know any simpler proof.

But if we multiply by another factor of $n!$,
then the result $(E_{2n}/n!)_{n \ge 0}$ {\em is}\/ a Hamburger
--- and indeed a Stieltjes --- moment sequence.
To see this, start by rewriting \reff{eq.star1} as
\be
   {E_{2n} \over (2n)!}
   \;=\;
   2 \sum\limits_{k=0}^\infty (-1)^k \, \bigl[ (k+\smhalf)\pi \bigr]^{-(2n+1)}
   \;.
 \label{eq.star1.bis}
\ee
Now multiply this by the Stieltjes integral representation
\be
   {(2n)! \over n!}
   \;=\;
   2^n \, (2n-1)!!
   \;=\;
   {2^n \over \sqrt{2\pi}}
   \int\limits_{-\infty}^\infty \! x^{2n} \, e^{-\half x^2} \, dx
\ee
to get
\be
   {E_{2n} \over n!}
   \;=\;
   {2^{n+1} \over \sqrt{2\pi}}
   \int\limits_{-\infty}^\infty \! dx \: e^{-\half x^2} \,
       \sum\limits_{k=0}^\infty {(-1)^k \over (k+\smhalf)\pi} \,
              \biggl( {x \over (k+\smhalf)\pi} \biggr) ^{\! 2n}
   \;.
\ee
Change variable to $y = x/[(k+\smhalf)\pi]$
and interchange integration and summation;  this leads to
\be
   {E_{2n} \over n!}
   \;=\;
   {2^{n+1} \over \sqrt{2\pi}}
   \int\limits_{-\infty}^\infty \! dy \: y^{2n} \:
       \sum\limits_{k=0}^\infty  (-1)^k \,
               \exp\bigl[- \smhalf (k+\smhalf)^2 \pi^2 y^2 \bigr]
   \;.
 \label{eq.density.thetalike}
\ee
The density here is positive because
each term with even $k$ dominates the term $k+1$:
\be
   \exp\bigl[- \smhalf (k+\smhalf)^2 \pi^2 y^2 \bigr]
   \;\ge\;
   \exp\bigl[- \smhalf (k+1+\smhalf)^2 \pi^2 y^2 \bigr]
   \;.
\ee
It follows that $(E_{2n}/n!)_{n \ge 0}$
is a Stieltjes moment sequence.\footnote{
   The fascinating article of Biane, Pitman and Yor \cite{Biane_01}
   discusses some densities that seem closely related to
   --- but apparently different from ---
   the one occurring in \reff{eq.density.thetalike}.
   I thank Christophe Vignat for drawing my attention to this article.
}
Its ordinary generating function is therefore given by a
Stieltjes-type continued fraction with coefficients $\alpha_i > 0$;
but no explicit formula for these coefficients seems to be known.

\section{Euler numbers, part 2}  \label{sec.euler2}

Thus far we have considered the sequence of
Euler numbers $E_n$ divided by factorials.
Now we consider the sequence of Euler numbers $E_n$ {\em tout court}\/,
along with its even and odd subsequences.

We have already seen that $(E_{n+1}/(n+1)!)_{n \ge 0}$
is a Hamburger moment sequence.
Since $((n+1)!)_{n \ge 0}$ is also a Hamburger (in fact a Stieltjes)
moment sequence,
it follows that their product $(E_{n+1})_{n \ge 0}$
is again a Hamburger moment sequence.
Similarly, we have seen that $(E_{2n+1}/(2n+1)!)_{n \ge 0}$
is a Stieltjes moment sequence;
and since ${((2n+1)!)_{n \ge 0}}$ is a Stieltjes moment sequence,
it follows that their product $(E_{2n+1})_{n \ge 0}$
is a Stieltjes moment sequence.
And finally, we have seen that $(E_{2n}/n!)_{n \ge 0}$
is a Stieltjes moment sequence;
and since $(n!)_{n \ge 0}$ is a Stieltjes moment sequence,
it follows that their product $(E_{2n})_{n \ge 0}$
is a Stieltjes moment sequence.
In this section we will obtain explicit expressions
for these sequences' representing measures
and for the continued-fraction expansions
of their ordinary generating functions.


Start by rewriting \reff{eq.star3b} as
\be
   {E_n \over n!}
   \;=\;
   2 \left[ \sum\limits_{k=0}^\infty
               \biggl( {2 \over (4k+1)\pi} \biggr) ^{\! n+1}
            \;-\;
            (-1)^n \sum\limits_{k=0}^\infty
               \biggl( {2 \over (4k+3)\pi} \biggr) ^{\! n+1}
     \right]
   \;.
\ee
Now multiply by the Stieltjes integral representation
$n! = \int\limits_0^\infty \! x^n \, e^{-x} \, dx$
to get
\begin{eqnarray}
   & &  \hspace*{-2cm}
   E_n
   \;=\;
   2 \left[ \int\limits_0^\infty \! dx \, e^{-x} \sum\limits_{k=0}^\infty
               {2 \over (4k+1)\pi}
               \biggl( {2x \over (4k+1)\pi} \biggr) ^{\! n}
     \right.
         \nonumber \\
   & & \qquad\qquad
     \left.
            \;-\; (-1)^n
            \int\limits_0^\infty \! dx \, e^{-x}
            \sum\limits_{k=0}^\infty
               {2 \over (4k+3)\pi}
               \biggl( {2x \over (4k+3)\pi} \biggr) ^{\! n}
     \right]
   \;.
\end{eqnarray}
Change variable to $y = 2x/[(4k+1)\pi]$ in the first term,
and $y = 2x/[(4k+3)\pi]$ in the second,
and interchange integration and summation;  this leads to
\begin{subeqnarray}
   E_n
   & = &
   2 \left[ \int\limits_0^\infty \! 
               {e^{-(\pi/2)y}  \over 1 - e^{-2 \pi y}}
               \: y^n \: dy
            \;-\;
            \int\limits_0^\infty \! 
               {e^{-(3\pi/2)y}  \over 1 - e^{-2 \pi y}}
               \: (-y)^n \: dy
     \right]
         \\[2mm]
   & = &
   \int\limits_{-\infty}^\infty \!
          y^n \: {e^{(\pi/2)y} \over \sinh \pi y} \: dy
   \;.
 \slabel{eq.En.hamburger0b}
 \label{eq.En.hamburger0}
\end{subeqnarray}
The integral \reff{eq.En.hamburger0b} is absolutely convergent for $n \ge 1$;
for $n=0$ it is valid as a principal-value integral at $y=0$.
In particular we have
\be
   E_{n+1}  \;=\;  \int\limits_{-\infty}^\infty \! 
          y^n \: {y \, e^{(\pi/2)y} \over \sinh \pi y} \: dy
   \;,
 \label{eq.En.hamburger}
\ee
which represents $E_{n+1}$ as the $n$th moment of a positive measure on $\R$.
Hence $(E_{n+1})_{n \ge 0}$ is a Hamburger moment sequence.
It is H-determinate by virtue of \reff{eq.euler.asymptotics}
and Corollary~\ref{cor.determinacy}.

Note also that multiplying \reff{eq.En.hamburger0b} by $t^n/n!$
and summing $\sum_{n=0}^\infty$, we recover the two-sided Laplace transform
\be
   \sec t + \tan t
   \;=\;
   \int\limits_{-\infty}^\infty \!
       {e^{(t + \pi/2)y} \over \sinh \pi y} \: dy
   \;,
\ee
which is valid for $-3\pi/2 < \real t < \pi/2$
as a principal-value integral \cite[6.2(8)]{Erdelyi_ET_54},
or equivalently (by symmetrizing)
\be
   \sec t + \tan t
   \;=\;
   \int\limits_{-\infty}^\infty \!
       {\sinh (t + \pi/2)y \over \sinh \pi y} \: dy
   \;.
\ee

For $n$ even we can combine the $y \ge 0$ and $y \le 0$ contributions
in \reff{eq.En.hamburger0b} to obtain
\cite[3.523.4]{Gradshteyn_8th_ed} \cite[24.7.6]{NIST}
\be
   E_{2n}  \;=\; \int\limits_0^\infty \!
          y^{2n} \: \sech\Bigl( {\pi \over 2} y \Bigr) \: dy
   \;,
 \label{eq.E2n.stieltjes}
\ee
while for $n$ odd a similar reformulation gives
\cite[3.523.2]{Gradshteyn_8th_ed}
\be
   E_{2n+1}  \;=\; \int\limits_0^\infty \!
          y^{2n} \; y \csch\Bigl( {\pi \over 2} y \Bigr) \: dy
   \;.
 \label{eq.E2n+1.stieltjes}
\ee
We have thus explicitly expressed $(E_{2n})_{n \ge 0}$
and $(E_{2n+1})_{n \ge 0}$ as Stieltjes moment sequences.
By \reff{eq.euler.asymptotics} and Theorem~\ref{thm.determinacy}(b)
they are S-determinate.
And since the measures in \reff{eq.E2n.stieltjes}/\reff{eq.E2n+1.stieltjes}
are continuous, Proposition~\ref{prop.determinacy.berg}
implies that these sequences are also H-determinate.

The moment representations \reff{eq.E2n.stieltjes}/\reff{eq.E2n+1.stieltjes}
can also be expressed nicely in terms of the
{\em Lerch transcendent}\/ (or {\em Lerch zeta function}\/)
\cite[\S9.55]{Gradshteyn_8th_ed} \cite[\S25.14]{NIST},
which we take to be defined by the integral representation
\be
   \Phi(z,s,\alpha)
   \;=\;
   {1 \over \Gamma(s)} \int\limits_0^\infty
       {t^{s-1} \, e^{-\alpha t} \over 1 \,-\, z e^{-t}}
       \: dt
 \label{def.lerch}
\ee
for $\real s > 0$, $\real \alpha > 0$, and $z \in \C \setminus [1,\infty)$.
For $|z| < 1$ we can expand the integrand in a Taylor series in $z$
and then interchange integration with summation:
this yields
\be
   \Phi(z,s,\alpha)
   \;=\;
   \sum_{n=0}^\infty {z^n \over (n+\alpha)^s}
   \;,
 \label{eq.lerch.series}
\ee
valid for $\real s > 0$, $\real \alpha > 0$, and $|z| < 1$.
Moreover, an application of Lebesgue's dominated convergence theorem
to the same series expansion shows that \reff{eq.lerch.series} holds
also for $|z| = 1$ with the exception of $z = 1$.\footnote{
   {\sc Proof.}
   $\displaystyle \left| \sum\limits_{n=0}^N u^n \right|
    \,=\, \left| {1 - u^{N+1} \over 1-u} \right|
    \,\le\, {2 \over |1-u|}$
   whenever $|u| \le 1$ and $u \neq 1$.
   Applying this with $u = z e^{-t}$ shows that
   the dominated convergence theorem applies to the
   Taylor expansion in $z$ whenever $|z| \le 1$ and $z \neq 1$.
}
And under the stronger hypothesis $\real s > 1$ we can take $z \uparrow 1$
and conclude that \reff{eq.lerch.series} holds also for $z = 1$.

Let us now use \reff{eq.lerch.series} for $z = \pm 1$:
then \reff{eq.star1.bis} and \reff{eq.star2} can be written as
\begin{eqnarray}
   {E_{2n} \over (2n)!}
   & = &
   {2 \over \pi^{2n+1}} \, \Phi(-1,2n+1,\smhalf)
        \\[2mm]
   {E_{2n+1} \over (2n+1)!}
   & = &
   {2 \over \pi^{2n+2}} \, \Phi(1,2n+2,\smhalf)
\end{eqnarray}
Using \reff{def.lerch} to express $\Gamma(s) \, \Phi(z,s,\alpha)$
as an integral, we recover
\reff{eq.E2n.stieltjes} and \reff{eq.E2n+1.stieltjes}.

\bigskip

We can also obtain continued fractions for the ordinary generating
functions of these three sequences.
For $(E_{n+1})_{n \ge 0}$ we have the Jacobi-type continued fraction
\be
   \sum_{n=0}^\infty E_{n+1} \, t^n
   \;=\;
   \cfrac{1}{1 - t - \cfrac{t^2}{1 - 2t - \cfrac{3 t^2}{1 - 3t -  \cfrac{6 t^2}{1- 4t - \cfrac{10 t^2}{1 - \cdots}}}}}
 \label{eq.En1.Jfrac}
\ee
with coefficients $\gamma_n = n+1$ and $\beta_n = n(n+1)/2$.
This continued fraction ought to be classical,
but the first mention of which I am aware
is a 2006 contribution to the OEIS by an amateur mathematician, Paul D.~Hanna,
who found it empirically \cite{Hanna_06};
it~was proven a few years later by Josuat-Verg\`es \cite{Josuat-Verges_14}
by a combinatorial method (which also yields a $q$-generalization).\footnote{
   {\bf Note Added (31~July 2018):}
   Jiang Zeng has informed me that this continued fraction
   {\em is}\/ classical!
   Stieltjes \cite[eq.~(14)]{Stieltjes_1890} showed in 1890 that
   the polynomials $Q_n(a,x)$ defined by the exponential generating function
   $\sum\limits_{n=0}^\infty Q_n(a,x) \, \displaystyle {t^n \over n!}
    = (\cos t - x \sin t)^{-a}$
   have the ordinary generating function
$$
   \sum\limits_{n=0}^\infty Q_n(a,x) \, t^n
   \;=\;
   \cfrac{1}{1 - axt - \cfrac{a(x^2+1)t^2}{1 - (a+2)xt - \cfrac{2(a+1)(x^2+1)t^2}{1 - (a+4)x - \cfrac{3(a+2)(x^2+1)t^2}{1 - \cdots}}}}
$$
   with coefficients $\gamma_n = (a+2n)x$ and $\beta_n = n(a+n-1)(x^2 + 1)$
   (see also \cite[Theorem~3.17]{Josuat-Verges_14}).
   Since
   $\sum\limits_{n=0}^\infty E_{n+1} \, \displaystyle {t^n \over n!}
    = (1 - \sin t)^{-1}$
   and $(\cos t - \sin t)^{-2} = (1 - \sin 2t)^{-1}$,
   we have $Q_n(2,1) = 2^n E_{n+1}$
   (see also \cite[Proposition~2.1]{Josuat-Verges_14}).
   So taking $a=2$ and $x=1$ in Stieltjes' formula
   and replacing $t$ by $t/2$ yields \reff{eq.En1.Jfrac}.
   Note also that taking $a=1$ and $x=1$ in Stieltjes' formula
   yields \reff{eq.SnJfrac} below;
   and taking $x=0$ in Stieltjes' formula
   yields \reff{def.secpow}/\reff{eq.secpow.contfrac} below.
     \label{footnote.Stieltjes_1890}
}

\bigskip

{\bf Remark.}
The J-fraction \reff{eq.En1.Jfrac} does not arise by contraction
from any S-fraction.
Indeed, if we use the contraction formula \reff{eq.contraction_even.coeffs}
and solve for $\balpha$,
we find $(\alpha_1,\alpha_2,\alpha_3,\alpha_4,\alpha_5)$ $= (1,1,1,3,0)$,
but then $\alpha_5 \alpha_6 = \beta_3 = 6$ has no solution.
\myendremark

\bigskip

For the even and odd subsequences, we have Stieltjes-type continued fractions:
\be
   \sum_{n=0}^\infty E_{2n} \, t^n
   \;=\;
   \cfrac{1}{1 - \cfrac{1^2 t}{1 - \cfrac{2^2 t}{1 -  \cfrac{3^2 t}{1- \cdots}}}}
 \label{eq.cfrac.secant}
\ee
with coefficients $\alpha_n = n^2$,
and
\be
   \sum_{n=0}^\infty E_{2n+1} \, t^n
   \;=\;
   \cfrac{1}{1 - \cfrac{1\cdot2 t}{1 - \cfrac{2\cdot3 t}{1 -  \cfrac{3\cdot4 t}{1- \cdots}}}}
 \label{eq.cfrac.tangent}
\ee
with coefficients $\alpha_n = n(n+1)$.
These formulae were found by Stieltjes \cite[p.~H9]{Stieltjes_1889} in 1889
and by Rogers \cite[p.~77]{Rogers_07} in 1907.
They were given beautiful combinatorial proofs by
Flajolet \cite{Flajolet_80} in 1980.


\bigskip

Since $(E_{n+1})_{n \ge 0}$ is a Hamburger moment sequence,
it is natural to ask about the full sequence $(E_n)_{n \ge 0}$.
Is it a Hamburger moment sequence?
The answer is no,
because the $3 \times 3$ Hankel matrix
$\displaystyle{ \begin{bmatrix}
                    E_0 & E_1 & E_2 \\
                    E_1 & E_2 & E_3 \\
                    E_2 & E_3 & E_4
                 \end{bmatrix}
              }
 =
 \displaystyle{ \begin{bmatrix}
                    1 & 1 & 1 \\
                    1 & 1 & 2 \\
                    1 & 2 & 5
                 \end{bmatrix}
              }$
has determinant $-1$.
But a much stronger result is true:

\begin{proposition}
   \label{prop.En.nonhamburger}
No arithmetic-progression subsequence $(E_{n_0 + jN})_{N \ge 0}$
with $n_0$ {\em even} and $j$ {\em odd}\/ is a Hamburger moment sequence.
\end{proposition}

\proof
For $n_0 \ge 1$, \reff{eq.En.hamburger} yields
\be
   E_{n_0+jN}  \;=\;  \int\limits_{-\infty}^\infty \!
          y^{jN} \: {y^{n_0} \, e^{(\pi/2)y} \over \sinh \pi y} \: dy
   \;.
\ee
When $n_0$ is even ($\ge 2$) and $j$ is odd,
this represents $(E_{n_0 + jN})_{N \ge 0}$
as the moments of a signed measure on $\R$ with exponential decay
that is {\em not}\/ a positive measure.
Corollary~\ref{cor.uniqueness} then implies that $(E_{n_0 + jN})_{N \ge 0}$
is not a Hamburger moment sequence.
The assertion for $n_0 = 0$ then follows from the assertion for $n_0 = 2j$.
\qed

Since $(E_{n+1})_{n \ge 0}$ is a Hamburger moment sequence
with a representing measure of infinite support,
it follows that all the Hankel determinants
$\Delta_n^{(m)} = \det (E_{i+j+m})_{0 \le i,j \le n-1}$
for $m$ {\em odd}\/ are strictly positive.
On the other hand, the $j=1$ case of Proposition~\ref{prop.En.nonhamburger}
implies that for every {\em even}\/ $m$
there must exist at least one $n$ such that $\Delta_n^{(m)} < 0$.
But which one(s)?
The question of the sign of $\Delta_n^{(m)}$ for $m$ even
seems to be quite delicate,
and I am unable to offer any plausible conjecture.
%

\bigskip

{\bf Remarks.}
1.  Although the sequence $(E_n)_{n \ge 0}$ of Euler numbers
is not a Stieltjes or even a Hamburger moment sequence,
it is log-convex.
This can be proven inductively from the recurrence \reff{eq.euler.recurrence}
\cite[Example~2.2]{Liu_07}.
Alternatively, it can be proven by observing that the
tridiagonal matrix \reff{eq.tridiag} associated to the
continued fraction \reff{eq.En1.Jfrac}
is totally positive of order 2,
i.e.\ $\beta_n \ge 0$, $\gamma_n \ge 0$
and $\gamma_n \gamma_{n+1} - \beta_{n+1} \ge 0$ for all $n$.
This implies \cite{Zhu_13,Sokal_totalpos} that
$(E_{n+1})_{n \ge 0}$ is log-convex.
And since $E_0 E_2 - E_1^2 = 0$, it follows that also $(E_n)_{n \ge 0}$
is log-convex.
See also \cite{Zhu_18a} for some stronger results.

2.  Dumont \cite[Proposition~5]{Dumont_95}
found a nice Jacobi-type continued fraction
also for the sequence of Euler numbers with some sign changes:
\be
   \sum_{n=0}^\infty (-1)^{n(n-1)/2} \, E_{n+1} \, t^n 
   \;=\;
   \cfrac{1}{1 - t + \cfrac{3t^2}{1 + \cfrac{5 t^2}{1 - t +  \cfrac{14 t^2}{1 + \cfrac{18t^2}{1 - t - \cdots}}}}}
 \label{eq.EnJfrac.2}
\ee
with coefficients $\gamma_{2k} = 1$, $\gamma_{2k+1} = 0$,
$\beta_{2k-1} = -k(4k-1)$ and $\beta_{2k} = -k(4k+1)$.
\myendremark

\section{Springer numbers}  \label{sec.springer}

We now turn to the sequence of Springer numbers.
Since $\cos t - \sin t = {\sqrt{2} \cos(t + \pi/4)}$,
the partial-fraction expansion \reff{eq.partfrac.sec} for secant yields
\be
   {1 \over \cos t - \sin t}
   \;=\;
   {1 \over \sqrt{2}}
  \lim\limits_{N \to \infty} \, \sum\limits_{k=-N}^N
        {(-1)^k \over (k+\smquarter)\pi - t}
   \;.
\ee
Inserting this into the exponential generating function \reff{def.springer}
of the Springer numbers
and extracting coefficients of powers of $t$ on both sides, we obtain
\be
   {S_n \over n!}
   \;=\;
   {1 \over \sqrt{2}}
   \sum\limits_{k=-\infty}^\infty
   (-1)^k \, \bigl[ (k+\smquarter)\pi \bigr]^{-(n+1)}
 \label{eq.springer.star1}
\ee
(with the interpretation $\lim\limits_{N \to \infty} \, \sum\limits_{k=-N}^N$
 when $n=0$).
Since \reff{eq.springer.star1} represents
every arithmetic-progression subsequence
$(\widetilde{S}_{n_0 + jN})_{N \ge 0}$ [where $\widetilde{S}_n = S_n/n!$]
as the moments of a signed measure supported on $[-4/3\pi, (4/\pi)^j]$
that is {\em not}\/ a positive measure,
it follows by Corollary~\ref{cor.uniqueness}
that no such sequence is a Hamburger moment sequence.

We now consider the sequence of Springer numbers {\em tout court}\/.
Start by rewriting \reff{eq.springer.star1} as
\be
   {S_n \over n!}
   \;=\;
   {1 \over \sqrt{2}}
     \left[ \sum\limits_{k=0}^\infty
                (-1)^k \, \bigl[ (k+\smquarter)\pi \bigr]^{-(n+1)}
            \;+\;
            (-1)^n \sum\limits_{k=0}^\infty
                (-1)^k \, \bigl[ (k+\smfrac{3}{4})\pi \bigr]^{-(n+1)}
     \right]
   \;.
\ee
Now multiply by the Stieltjes integral representation
$n! = \int\limits_0^\infty \! x^n \, e^{-x} \, dx$
to get
\begin{eqnarray}
   & &  \hspace*{-2cm}
   S_n
   \;=\;
   {1 \over \sqrt{2}}
   \left[ \int\limits_0^\infty \! dx \, e^{-x} \sum\limits_{k=0}^\infty
               {(-1)^k \over (k+\smquarter)\pi}
               \biggl( {x \over (k+\smquarter)\pi} \biggr) ^{\! n}
     \right.
         \nonumber \\
   & & \qquad\qquad
     \left.
            \;+\; (-1)^n
            \int\limits_0^\infty \! dx \, e^{-x}
            \sum\limits_{k=0}^\infty
               {(-1)^k \over (k+\smfrac{3}{4})\pi}
               \biggl( {x \over (k+\smfrac{3}{4})\pi} \biggr) ^{\! n}
     \right]
   \;.
\end{eqnarray}
Change variable to $y = x/[(k+\smquarter)\pi]$ in the first term,
and $y = x/[(k+\smfrac{3}{4})\pi]$ in the second,
and interchange integration and summation;  this leads to
\begin{subeqnarray}
   S_n
   & = &
   {1 \over \sqrt{2}}
   \left[ \int\limits_0^\infty \! 
               {e^{-(\pi/4)y}  \over 1 + e^{-\pi y}}
               \: y^n \: dy
            \;+\;
            \int\limits_0^\infty \! 
               {e^{-(3\pi/4)y}  \over 1 + e^{-\pi y}}
               \: (-y)^n \: dy
     \right]
         \\[2mm]
   & = &
   {1 \over 2 \sqrt{2}}
   \int\limits_{-\infty}^\infty \!
          y^n \: {e^{(\pi/4)y} \over \cosh (\pi y/2)} \: dy
   \;,
 \slabel{eq.Sn.hamburger0b}
 \label{eq.Sn.hamburger0}
\end{subeqnarray}
which is absolutely convergent for all $n \ge 0$.
It follows that $(S_n)_{n \ge 0}$ is a Hamburger moment sequence.
It is H-determinate by virtue of \reff{eq.springer.asymptotics}
and Corollary~\ref{cor.determinacy}.
Since the unique representing measure has support equal to all of $\R$,
it follows that $(S_n)_{n \ge 0}$ is not a Stieltjes moment sequence.
This can alternatively be seen from the fact that
the $3 \times 3$ once-shifted Hankel matrix
$\displaystyle{ \begin{bmatrix}
                    S_1 & S_2 & S_3 \\
                    S_2 & S_3 & S_4 \\
                    S_3 & S_4 & S_5
                 \end{bmatrix}
              }
 =
 \displaystyle{ \begin{bmatrix}
                    1 & 3 & 11 \\
                    3 & 11 & 57 \\
                    11 & 57 & 361
                 \end{bmatrix}
              }$
has determinant~$-96$.

Note also that multiplying \reff{eq.Sn.hamburger0b} by $t^n/n!$
and summing $\sum_{n=0}^\infty$, we recover the two-sided Laplace transform
\be
   {1 \over \cos t - \sin t}
   \;=\;
   {1 \over 2 \sqrt{2}}
   \int\limits_{-\infty}^\infty \!
       {e^{(t + \pi/4)y} \over \cosh(\pi y/2)} \: dy
   \;,
\ee
which is valid for $-3\pi/4 < \real t < \pi/4$ \cite[6.2(11)]{Erdelyi_ET_54}.

For $n$ even we can combine the $y \ge 0$ and $y \le 0$ contributions
in \reff{eq.Sn.hamburger0b}
to obtain
\be
   S_{2n}
   \;=\;
   {1 \over \sqrt{2}}
   \int\limits_0^\infty \!
          y^{2n} \: {\cosh (\pi y/4)
                     \over
                     \cosh (\pi y/2)
                    } \: dy
   \;,
 \label{eq.S2n.stieltjes}
\ee
while for $n$ odd a similar reformulation gives
\be
   S_{2n+1}  \;=\;
   {1 \over \sqrt{2}}
   \int\limits_0^\infty \!
          y^{2n} \: {y \, \sinh (\pi y/4)
                     \over
                     \cosh (\pi y/2)
                    } \: dy
   \;.
 \label{eq.S2n+1.stieltjes}
\ee
We have thus explicitly expressed $(S_{2n})_{n \ge 0}$
and $(S_{2n+1})_{n \ge 0}$ as Stieltjes moment sequences.
By \reff{eq.springer.asymptotics} and Theorem~\ref{thm.determinacy}(b)
they are S-determinate.
And since the measures in \reff{eq.S2n.stieltjes}/\reff{eq.S2n+1.stieltjes}
are continuous, Proposition~\ref{prop.determinacy.berg}
implies that these sequences are also H-determinate.

We can also obtain continued fractions for the ordinary generating
functions of these three sequences.
For $(S_n)_{n \ge 0}$ we have the Jacobi-type continued fraction
\be
   \sum_{n=0}^\infty S_n \, t^n
   \;=\;
   \cfrac{1}{1 - t - \cfrac{2\cdot 1^2 \, t^2}{1 - 3t - \cfrac{2\cdot 2^2 \, t^2}{1 - 5t -  \cfrac{2 \cdot 3^2 \, t^2}{1- 7t - \cfrac{2 \cdot 4^2 \, t^2}{1 - \cdots}}}}}
 \label{eq.SnJfrac}
\ee
with coefficients $\gamma_n = 2n+1$ and $\beta_n = 2n^2$.
This formula was proven a few years ago
by Josuat-Verg\`es \cite{Josuat-Verges_14},
by a combinatorial method that also yields a $q$-generalization;
it~was independently found (empirically) by an amateur mathematician,
Sergei N.~Gladkovskii \cite{Gladkovskii_13}.\footnote{
   {\bf Note Added (31~July 2018):}
   This formula is also a special case of a result
   of Stieltjes \cite{Stieltjes_1890}:
   see footnote~\ref{footnote.Stieltjes_1890} above.
}
The fact that $\beta_n > 0$ for all $n$ tells us again that
$(S_n)_{n \ge 0}$ is a Hamburger moment sequence.

For the even Springer numbers we have the Stieltjes-type continued fraction
\cite[Corollary~3.3]{Dumont_95}
\be
   \sum_{n=0}^\infty S_{2n} t^n
   \;=\;
   \cfrac{1}{1 - \cfrac{1\cdot3 t}{1 - \cfrac{4\cdot 4t}{1 - \cfrac{5\cdot 7 t}{1- \cfrac{8\cdot8 t}{1-\cdots}}}}}
\ee
with coefficients $\alpha_{2k-1} = (4k-3)(4k-1)$ and $\alpha_{2k} = (4k)^2$.
For the odd Springer numbers we have the Jacobi-type continued fraction
\cite{Gladkovskii_12}
\be
   \sum_{n=0}^\infty S_{2n+1} t^n
   \;=\;
   \cfrac{1}{1 - 11t - \cfrac{16\cdot1\cdot3\cdot5t^2}{1 - 75t - \cfrac{16\cdot4\cdot7\cdot9 t^2}{1 - 203t - \cfrac{16\cdot9\cdot11\cdot13 t}{1- \cdots}}}}
   \label{eq.springerodd.Jfrac}
\ee
with coefficients 
$\gamma_n  = 32n^2 + 32n + 11$ and $\beta_n = (4n-1)(4n)^2(4n+1)$.
This formula can be obtained as a specialization of a result
of Stieltjes \cite{Stieltjes_1890} (see \cite{Bala_17});
it can alternatively be obtained from
\cite[Propositions~7 and 8]{Dumont_95}
by the transformation formula for Jacobi-type continued fractions
under the binomial transform
\cite[Proposition~4]{Barry_09} \cite{Sokal_totalpos}.
Since the odd Springer numbers are a Stieltjes moment sequence,
their ordinary generating function is also given by a
Stieltjes-type continued fraction with coefficients $\alpha_i > 0$;
these coefficients can in principle be obtained from
\reff{eq.springerodd.Jfrac} by solving \reff{eq.contraction_even.coeffs},
but no explicit formula for them seems to be known
(and maybe no simple formula exists).

\bigskip

{\bf Remarks.}
1.  Although the sequence $(S_n)_{n \ge 0}$ of Springer numbers
is not a Stieltjes moment sequence, it is log-convex.
This follows \cite{Zhu_13,Sokal_totalpos}
from the fact that the
tridiagonal matrix \reff{eq.tridiag} associated to the
continued fraction \reff{eq.SnJfrac} is totally positive of order 2.
The log-convexity of the Springer numbers was conjectured a few years ago
by Sun \cite[Conjecture~3.4]{Sun_13}
and proven recently by Zhu {\em et al.}\/ \cite{Zhu_18b}
as a special case of a more general result.

2. Dumont \cite[Corollary~3.2]{Dumont_95}
also found a nice Jacobi-type continued fraction
for the sequence of Springer numbers with some sign changes:
\be
   \sum_{n=0}^\infty (-1)^{n(n-1)/2} \, S_n \, t^n 
   \;=\;
   \cfrac{1}{1 - t + \cfrac{4t^2}{1 - t + \cfrac{16 t^2}{1 - t +  \cfrac{36 t^2}{1- \cdots}}}}
 \label{eq.SnJfrac.2}
\ee
with coefficients $\gamma_n = 1$ and $\beta_n = -4n^2$.
This formula follows from \reff{eq.cfrac.secant} with $t$ replaced by $4t^2$,
combined with the identity
\be
   (-1)^{n(n-1)/2} \, S_n
   \;=\;
   \sum_{k=0}^{\lfloor n/2 \rfloor}
   \binom{n}{2k} \, (-4)^{k} \, E_{2k}
\ee
(which follows from the exponential generating functions
 \cite[p.~275]{Dumont_95})
and a general result about how Jacobi-type continued fractions
behave under the binomial transform
\cite[Proposition~4]{Barry_09} \cite{Sokal_totalpos}.
\myendremark

\section{What next?}

Enumerative combinatorialists are not content with merely counting sets;
we want to refine the counting by measuring one or more statistics.
To take a trivial example, an $n$-element set has $2^n$ subsets,
but we can classify these subsets according to their cardinality,
and say that there are $\binom{n}{k}$ subsets of cardinality $k$.
We then collect these refined counts in a {\em generating polynomial}\/
\be
   P_n(x)  \;\eqdef\; \sum_{A \subseteq [n]} x^{|A|}
           \;=\;  \sum_{k=0}^n \binom{n}{k} x^k
   \;,
\ee
which in this case of course equals $(1+x)^n$.
To take a less trivial example, the number of ways of partitioning
an $n$-element set into nonempty blocks is given by the Bell number $B_n$;
but we can refine this classification by saying that the
number of ways of partitioning
an $n$-element set into $k$ nonempty blocks is given by the Stirling number
$\stirlingsubset{n}{k}$, and then form the Bell polynomial
\be
   B_n(x)  \;=\;  \sum_{k=0}^n \stirlingsubset{n}{k} x^k
   \;.
\ee
We can then study generating functions, continued-fraction expansions,
moment representations and so forth for $B_n(x)$,
generalizing the corresponding results for $B_n = B_n(1)$.

In a similar way, the Euler and Springer numbers can be refined into polynomials
that count alternating permutations or snakes of type $B_n$
according to one or more statistics.
For instance, consider the polynomials $E_{2n}(x)$ defined by
\be
   (\sec t)^x  \;=\;  \sum_{n=0}^\infty E_{2n}(x) \, {t^{2n} \over (2n)!}
 \label{def.secpow}
\ee
where $x$ is an indeterminate.
They satisfy the recurrence \cite[p.~123]{Johnson_00}
\be
   E_{2n+2}(x)
   \;=\;
   x \sum_{k=0}^n \binom{2n+1}{2k} \, E_{2n-2k-1} \, E_{2k}(x)
\ee
with initial condition $E_0(x) = 1$.
It follows that $E_{2n}(x)$ is a polynomial of degree $n$ with
nonnegative integer coefficients,
which we call the \emph{secant power polynomial}.
The first few secant power polynomials are
\cite[A088874/A085734/A098906]{OEIS}
\begin{subeqnarray}
   E_0(x)  & = &   1  \\[1mm]
   E_2(x)  & = &   \phantom{1 + } x  \\[1mm]
   E_4(x)  & = &   \phantom{1 + } 2x + 3x^2  \\[1mm]
   E_6(x)  & = &   \phantom{1 + } 16 x + 30 x^2 + 15 x^3  \\[1mm]
   E_8(x)  & = &   \phantom{1 + } 272 x + 588 x^2 + 420 x^3 + 105 x^4
\end{subeqnarray}
Since $E_{2n}(1) = E_{2n}$,
these are a polynomial refinement of the secant numbers.
Moreover, since $\tan' = \sec^2$ and $(\log\sec)' = \tan$,
we have $E_{2n}(2) = E_{2n+1}$ and $E'_{2n}(0) = E_{2n-1}$,
so these are also a polynomial refinement of the tangent numbers.
Carlitz and Scoville \cite{Carlitz_75} proved that $E_{2n}(x)$
enumerates the alternating (down-up) permutations of $[2n]$ or $[2n+1]$
according to the number of records:
\be
   \sum_{\sigma \in \Alt_{2n}} \!\! x^{\rec(\sigma)}
   \;=\;
   E_{2n}(x)
 \label{eq.carlitz-scoville.1}
\ee
and
\be
   \sum_{\sigma \in \Alt_{2n+1}} \!\!\!\! x^{\rec(\sigma)}
   \;=\;
   x \, E_{2n}(1+x)
   \;.
 \label{eq.carlitz-scoville.2}
\ee
Here a {\em record}\/ (or {\em left-to-right maximum}\/)
of a permutation $\sigma \in \Sym_n$
is an index $i$ such that $\sigma_j < \sigma_i$ for all $j < i$.
(In particular, when $n \ge 1$, the index 1 is always a record.
 This explains why the $E_{2n}(x)$ for $n > 0$ start at order $x$.)

It turns out that
the ordinary generating function of the secant power polynomials
is given by a beautiful Stieltjes-type continued fraction,
which was found more than a century ago
by Stieltjes \cite[p.~H9]{Stieltjes_1889} and Rogers \cite[p.~82]{Rogers_07}
(see also \cite{Barry_12b,Josuat-Verges_14}):
\be
   \sum_{n=0}^\infty E_{2n}(x) \, t^n
   \;=\;
   \cfrac{1}{1 - \cfrac{x t}{1 - \cfrac{2(x+1) t}{1 -  \cfrac{3(x+2) t}{1- \cdots}}}}
 \label{eq.secpow.contfrac}
\ee
with coefficients $\alpha_n = n(x+n-1)$.
When $x=1$ this reduces to the expansion \reff{eq.cfrac.secant}
for the secant numbers;
when $x=2$ it becomes the expansion \reff{eq.cfrac.tangent}
for the tangent numbers.

The nonnegativity of the coefficients $\alpha_n$ in \reff{eq.secpow.contfrac}
for $x \ge 0$ implies, by Theorem~\ref{thm.stieltjes}, that for every $x \ge 0$,
the sequence $(E_{2n}(x))_{n \ge 0}$ is a Stieltjes moment sequence.
In~fact, $(E_{2n}(x))_{n \ge 0}$ has the explicit
Stieltjes moment representation \cite[pp.~179--181]{Chihara_78}\footnote{
   See \cite[eq.~5.13.2]{NIST} for the normalization.
}
\be
   E_{2n}(x)
   \;=\;
   {2^{x-1} \over \pi \, \Gamma(x)}
   \int\limits_{0}^\infty s^{2n} \:
          \biggl| \Gamma \Bigl( \displaystyle {x+is \over 2} \Bigr) \biggr|^2
         \: ds
   \;,
\ee
which reduces to \reff{eq.E2n.stieltjes} when $x=1$,
and to \reff{eq.E2n+1.stieltjes} when $x=2$.

The continued fraction \reff{eq.secpow.contfrac}
also implies, by Theorem~\ref{thm.stieltjes},
that for every $x \ge 0$,
every minor of the Hankel matrix $(E_{2i+2j}(x))_{i,j \ge 0}$
is a nonnegative real number.
But a vastly stronger result turns out to be true \cite{Sokal_totalpos}:
namely, every minor of the Hankel matrix $(E_{2i+2j}(x))_{i,j \ge 0}$
is a polynomial in $x$ with nonnegative integer coefficients!
This {\em coefficientwise}\/ Hankel-total positivity
arises in a wide variety of sequences of combinatorial polynomials
(sometimes in many variables)
--- in some cases provably, in other cases conjecturally.
But that is a story for another day.


\section*{Acknowledgments}

I wish to thank Christian Berg, Gerald Edgar, Sergei Gladkovskii, Paul Hanna,
Matthieu Josuat-Verg\`es, Mathias P\'etr\'eolle,
Christophe Vignat, Jiang Zeng and Bao-Xuan Zhu
for helpful conversations and/or correspondence.
I am also grateful to an anonymous referee for several valuable suggestions.

This research was supported in part by
Engineering and Physical Sciences Research Council grant EP/N025636/1.

\end{document}